\documentclass{svjour3}                     % onecolumn (standard format)
\smartqed  % flush right qed marks, e.g. at end of proof
\usepackage{graphicx}
\usepackage{mathptmx}      % use Times fonts if available on your TeX system
\usepackage[T1]{fontenc}
\usepackage[latin1]{inputenc}
\usepackage{amsmath}
\usepackage{amsfonts}
\usepackage{amssymb}
\usepackage{graphicx,color}
\usepackage{float}
\usepackage{geometry}
\usepackage{subfig}
\usepackage[numbers]{natbib}
\newcommand\T{\rule{0pt}{2.6ex}}

\newcommand\B{\rule[-1.2ex]{0pt}{0pt}}
\usepackage{pslatex}

\newcommand{\Rn}[1][m]{\ensuremath{\mathbf{R}^{#1}}}
\DeclareMathOperator{\diag}{diag}
\journalname{Computaional Geosciences}
\begin{document}

\title{A global method for coupling transport
 with chemistry in heterogeneous porous media%\thanks{Grants or other notes
}


\author{Laila AMIR  \and Michel KERN %etc.
}


\institute{ 
Laila AMIR \and Michel KERN \at
	      INRIA-Rocquencourt, B.P. 105, F-78153 Le Chesnay Cedex, France.\\
              \email{Michel.kern@inria.fr} 
\and Laila AMIR\at
              ITASCA Consultants, S.A, 64, Chemin des Mouilles 69134 -
              Ecully, France, \\
              \email{Laila.Amir@inria.fr}
\and   
}

\date{Received: date / Accepted: date}

\maketitle

\begin{abstract}
  Modeling reactive transport in porous media, using a local chemical
  equilibrium assumption, leads to a system of advection-diffusion
  PDE's coupled with algebraic equations.  When solving this coupled
  system, the algebraic equations have to be solved at each grid point
  for each chemical species and at each time step.  This leads to a
  coupled non-linear system.  In this paper a global solution approach
  that enables to keep the software codes for transport and chemistry
  distinct is proposed. The method applies the Newton-Krylov framework
  to the formulation for reactive transport used in operator
  splitting.  The method is formulated in terms of total mobile and
  total fixed concentrations and uses the chemical solver as a black
  box, as it only requires that on be able to solve chemical
  equilibrium problems (and compute derivatives), without having to
  know the solution method. An additional advantage of the
  Newton-Krylov method is that the Jacobian is only needed as an
  operator in a Jacobian matrix times vector product. The proposed
  method is tested on the MoMaS reactive transport benchmark.

\keywords{Geochemistry \and transport in porous media \and
  Newton-Krylov methods \and advection--diffusion--reaction equations}
\subclass{76V05 \and 65M99}
\end{abstract}

\section{Introduction}
\label{intro}
The simulation of multi-species reacting systems in porous media is of
importance in several different fields: for computing the near field
in nuclear waste simulations, in the treatment of bio-remediation, in
CO2 sequestration simulations and in the evaluation of underground
water quality.

This work deals with numerical methods for solving coupled transport
and chemistry problems. The transport of solutes in porous media is
described by partial differential equations of advection--diffusion
type, wheres multi-species chemistry involves the solution of
ordinary differential equations (if the reactions are kinetic) or
nonlinear algebraic equations (if local equilibrium is assumed). After
discretization, one is led to a system of nonlinear equations, coupled
the unknowns for all chemical species at all grid points. 

After the influential paper by Yeh and Tripathi~\cite{yehtrip},
operator splitting methods, where transport and chemistry are solved
for separately at each time step (possibly iterating to convergence),
became the methods of choice. Some representative papers where
operator splitting methods are used are
~\cite{sacaay:01},~\cite{sacaay:00},~\cite{carrmosbeh:04}, \cite{kamike:03},
~\cite{LuBuOl:00},~\cite{samxuyang:08}.  Operator splitting methods
are easy to implement, and the splitting errors can be controlled by
carefully restricting the time step. On the other hand, the time-step
restriction can become their main drawback, as it can be difficult to
get the fixed point iteration to converge for more difficult problems.

More recently, global methods have become popular, due to the increase
in computing power now available. In this approach, the full
non-linear system is solved in one step, usually by some form of
Newton's method. Most papers use the Direct Substitution Approach (see
~\cite{hammvallicht05},~\cite{FahCarYouAck:08}), where one
\emph{substitutes} the chemical equations in the transport equations.
On the other hand, the problem can also be put in the form of a
Differential Algebraic Equations (DAE), enabling the use of powerful
software (see~\cite{dieulerhkern:09}). Finally, the chemical equations
can be eliminated locally, and a system involving transport equations,
with a source term coming from the reactions has to be solved. This
approach is taken in~\cite{knabner05,knabner07}, where additionally a
reduction method leads to a smaller system.  Most of the papers quoted
above employ a Newton method for solving the nonlinear system at each
time step, with the difficulty that the Jacobian matrix has to be
computed, stored and factored. This can become problematic for large
problems, and Hammond et al.~\cite{hammvallicht05} have used the
Jacobian-Free Newton--Krylov method, where the Newton correction is
solved for by an iterative method. The Jacobian is only needed through
the computation of a directional derivative. The method keeps the fast
convergence of Newton's method, while only requiring Jacobian
matrix--vector products, and these can be approximated by finite
differences.

The method presented in this paper is a global method where the
chemical equations are eliminated locally, leading to a nonlinear
system where the transport and chemistry subsystems remain separated.
Thus the residual can be evaluated by calling separately written
transport and chemistry modules.  The system is then solved by a
Newton-Krylov method, and it will be shown how the Jacobian
matrix--vector product can also be computed by the same module. Thus
the main contribution of this paper is to show that a global method
can be implemented while still keeping transport and chemistry modules
separated. This property will be referred to as using ``black--box
solvers''. As the chemical equilibrium equations are not substituted
in the transport equations, the transport and chemistry parts of the
nonlinear residual are easily identified, and can each be computed by
calling on standard solution modules.

\medskip
An outline of the paper is as follows. In section~\ref{sec:TReqs} the
chosen model is explained, and the methods used for solving the
(non-reactive) transport part, and the chemical equilibrium
system are detailed. Section~\ref{sec:couplmodel} shows how we obtain the coupled
model. Couple formulations and coupling algorithms are the subject of
section~\ref{sec:couplalg}, beginning with a review of existing
methods, while our approach is presented in
section~\ref{sec:NewtKryl}. Numerical results, in particular
experience with the MoMaS benchmark, are shown in section~\ref{sec:numres}.

\section{Reactive transport equations}
\label{sec:TReqs}

In this work, the transport of several reacting species in a single
phase flow through a porous medium is considered. The species can
react both between themselves and with the porous matrix. In this
section the numerical methods used to
solve the individual subsystems of the coupled problem will be described.

\subsection{Transport model}
The transport of a single species through a porous medium (a domain
$\Omega \subset \mathbf{R}^d$, with $d=1,2$ or 3), with porosity
$\phi$, in a known Darcy field $u$, subject to dispersion and
molecular diffusion, follows the linear advection--dispersion equation
\begin{equation}
\label{eq:transp}
   \phi\,\dfrac{\partial c}{\partial t} + L(c) = q, \quad \text{in } \Omega
\end{equation}
where 
$$
L(c)=\nabla \cdot (u\, c)-\nabla\cdot(D \nabla c), 
$$
is the transport operator, and $q$ is a source term. The
diffusion--dispersion tensor $D$ is given by
\begin{equation*}
  D = d_e I + |u| \left(\alpha_L E(u) + \alpha_t (I-E(u))  \right),
  \quad E_{ij}(u) = \dfrac{u_i u_j}{|u|^2}, 
\end{equation*}
where $d_e$ is the molecular diffusion coefficient, and $\alpha_L$
(resp. $\alpha_t$) is the longitudinal (resp. transverse) dispersivity
coefficient.

In this work, we restrict to a one dimensional problem, so that the
transport equation over a bounded interval $\Omega=]0,L[$  can be
written as 
\begin{equation}
\begin{array}{l}
\phi\dfrac{\partial c}{\partial t}+\dfrac{\partial}{\partial x}
 \left( -D\dfrac{\partial c}{\partial x}+uc \right)= q, \quad 0 < x <
 L, \quad  0 < t < T,
\end{array}
\label{eq1}
\end{equation}
where the porosity $\phi$ and the diffusion--dispersion coefficient
$D$ can both depend on space. Because the flow is assumed
compressible, the velocity $u$ is taken to be a constant. 

The initial condition is $c(x,0)=c_0(x)$ and, in view of the
applications, the boundary conditions are a Dirichlet condition (given concentration)
$c(0,t)=c_d(t)$ at the left boundary  ($x=0$) and zero diffusive flux
$\dfrac{\partial c}{\partial x} = 0$ at the right
boundary ($x=L$). More general boundary conditions could easily be accommodated.

\subsubsection{Discretization in space}
We treat the space and time discretization separately, as we will use
different time discretizations for the different parts of the
transport operator. 

For space discretization a cell-centered finite volume scheme will be
used, see for instance~\cite{eymagallherb:00}.  The interval $[0,L]$
is divided into $N_g$ intervals
$[x_{i-\frac{1}{2}},x_{i+\frac{1}{2}}]$ of length $h_i$, where
$x_{\frac{1}{2}}=0,x_{N_g+\frac{1}{2}}=L$.  For $i=1,..,N_g$, denote
by $x_i$ the center and $x_{i+1/2}$ the right end of element $i$.
Finally, denote by $c_i$, $i=1,..,N_g$ the approximate solution in
cell $i$.

Equation~\eqref{eq1} is written in the form
\begin{equation}
  \label{eq:tr1dflux}
  \phi \dfrac{\partial c}{\partial t} + \dfrac{\partial
    \varphi}{\partial x} = q, 
\end{equation}
where the flux $\varphi(x,t)=-D\dfrac{\partial c}{\partial x}+uc$ has
been split as the sum of a diffusive flux $\varphi_d=-D\dfrac{\partial
  c}{\partial x}$ and an advective flux $\varphi_a=uc$.

Equation~\eqref{eq:tr1dflux} is integrated over a cell $]x_{i-1/2},
x_{i+1/2}[$ giving
\begin{equation}
\label{eq:trdisc}
\phi_i h_i \dfrac{d c_i}{dt} + \varphi_{d,i+\frac{1}{2}} +
\varphi_{a,i+\frac{1}{2}} - \varphi_{d,i-\frac{1}{2}} -
\varphi_{a,i-\frac{1}{2}} = h_i q_i, \qquad i=2, \dots, N_g. 
\end{equation}
The flux approximations required to close the system are provided by
finite differences. The diffusive flux needs a value for the diffusion
coefficient, which is taken as the harmonic average
(as done in mixed finite element methods):
\begin{equation}
\label{eq:fluxdisc}
\varphi_{d,i+\frac{1}{2}}=-D_{i+\frac{1}{2}}\left(\frac{c_{i+1} -
    c_i}{h_{i+\frac{1}{2}}}\right)   
\end{equation}
with
$$
D_{i+\frac{1}{2}}=\frac{2D_iD_{i+1}}{D_i+D_{i+1}}, \quad
D_{\frac{1}{2}}=D_1, \quad D_{N_g+\frac{1}{2}}=D_{N_g}
 \quad \mbox{and}\quad  h_{i+\frac{1}{2}}=\frac{h_i+h_{i+1}}{2}
$$ 
For the advective flux, an upwind approximation is used, so that
(assuming $u>0$), $\varphi_{a,i+\frac{1}{2}}=uc_i$

These approximations are corrected to take into account the boundary
conditions, both at $x=0$ and at $x=L$. The semi-discrete system can
be summarized by the finite dimensional system
\begin{equation}
  \label{eq:trans_ode}
  M \dfrac{dc}{dt} + L c = q + g, 
\end{equation}
where $c \in \Rn[N_g]$ now represents the vector of cell
concentrations, $L \in \Rn[N_g, N_g]$ is the matrix form of the
transport operator, $M \in \Rn[N_g, N_g]$ is a mass matrix accounting
for variable porosity and mesh size, $q \in \Rn[N_g]$ is a give source
term and $g \in \Rn[N_g]$ represents the effects of the boundary conditions.

\subsubsection{Time discretization}

Let denote by $\Delta t$ the time step (taken constant for simplicity)
used to discretize the time interval $[0, T]$ and denote by $c_i^n$
the (approximate) value of $c_i(n \Delta t)$.  The first and most
straightforward alternative is to discretize
equation~(\ref{eq:trans_ode}) by the backward Euler method, see for
instance~\cite{acher:NMEDE}. This is the method that is used in
section 3 to keep the description simple, but is not the recommended
method, as it leads to an overly diffusive scheme.

Better alternatives are obtained by exploiting the structure of the
transport operator, and by using different time discretizations for
the advective and for the diffusive parts. Specifically, the
diffusive terms should be treated implicitly, and the advective terms
are better handled explicitly.

If this idea is applied directly to equation~(\ref{eq:trans_ode}), the
resulting fully discrete scheme is only stable under a CFL
(Courant--Friedrichs--Lewy) condition $u \Delta t \leq \max_I h_i$. 
As this may be too severe a restriction (some of our applications
require integration over a very large time interval), an alternative
is to use an operator splitting scheme, as proposed by Siegel et
al.~\cite{siegmosackjaff:97} (see
also~\cite{hotackmos:04,MazBerPut:00}). In this work, splitting is
used only within the (linear) transport step, but recent papers by
Ka\v{c}ur et al.~\cite{KacFrol:06,KacMalRem:05} apply splitting
directly to a transport with sorption model by solving (analytically)
a nonlinear advection step, followed by a nonlinear diffusion step.
This is different from operator splitting as used in geochemical
models, as the chemistry terms are solved for together with the
transport terms.

The splitting scheme works by taking several small time steps of
advection, controlled by a CFL condition, within a large time step of
diffusion. The scheme has been shown to be unconditionally stable, and
has a good behavior in advection dominated situations.

More precisely, the time step $\Delta t$ will be used as the diffusion
time step, it is divided into M time steps of advection $\Delta t_c$
such that $\Delta t=M\Delta t_c$ where M >1, the advection time step
will be controlled by CFL condition. Equation~\eqref{eq:tr1dflux} will
be solved over the time step $[t^n, t^{n+1}]$ by first solving the
advection equation $\phi\dfrac{\partial c}{\partial
  t}+\dfrac{\partial}{\partial x} (uc)=0$ over $M$ steps of size
$\Delta t_c$ each, and then solving the diffusion equation
$\;\phi\dfrac{\partial c}{\partial t}+\dfrac{\partial}{\partial x}
(-D\dfrac{\partial c}{\partial x})=q$ starting from the value at the
end of the advection step.

\paragraph{Advection step}

The interval $[t^{n},t^{n+1}]$ is divided into $M$ intervals
$[t^{n,m},t^{n,m+1}]$, $m=0,...M-1$, where  
$t^{n,0}=t^{n},\,t^{n,M}=t^{n+1}$. Denote $c_{i}^{n,m}$ the
approximate concentration c at time $t^{n,m}$ and $c^{n,0}=c^{n}$. The
advection equation is discretized in time using the explicit Euler
method to obtain  
\begin{equation}
\left\{
  \begin{aligned}
     &\phi_{i}(\frac{c_{i}^{n,m+1}-c_{i}^{n,m}}{\Delta t_c})+ u 
    \left( \dfrac{c_{i}^{n,m} - c_{i-1}^{n,m}}{h_{i-1/2}}\right) =0,
    \quad i=2, \dots, N_g, \\
    & c_1^{n, m+1} = c_g(t^{n, m+1})
  \end{aligned}
\right.
 \quad m=0, \dots, M-1.
\end{equation}

\paragraph{Diffusion step}
The diffusion part is discretized by an implicit Euler scheme,
starting from $c_i^{n, M}$ :
\begin{equation}
-\frac{D_{i-\frac{1}{2}}}{h_{i-\frac{1}{2}}}{\Delta t} c_{i-1}^{n+1}
+ \left(\phi_i{h_i}+\frac{D_{i+\frac{1}{2}}}{h_{i+\frac{1}{2}}}{\Delta t} +
\frac{D_{i-\frac{1}{2}}}{h_{i-\frac{1}{2}}}{\Delta t}\right) c_{i}^{n+1} - 
\frac{D_{i+\frac{1}{2}}}{h_{i+\frac{1}{2}}}{\Delta t}c_{i+1}^{n+1} = 
\phi_i{h_i}c_{i}^{n,M}+q_ih_i\Delta t, \quad i=2, \dots, N_g-1 
\end{equation}
As above, 2 equations accounting for the boundary conditions must be added.

\subsection{Chemical equations}
\label{sec:chimie}

The chemical model is described in this section. In this study, we
assume a local chemical equilibrium at every point, which means that
the chemical phenomena occur on much faster scale than transport
phenomena. This is a common modeling assumption for reactive
transport in porous media, at least when the only reactions
considered are aqueous phase and sorption reactions (these are
``sufficiently fast'' reactions according to
Rubin~\cite{rubin:83}). This would not be the case if mineral
dissolution was taken into account, as these reactions typically need
kinetic models. 

Consider a set of $N_e$ chemical species $(X_j)_{j=1,\dots, N_e}$
linked by $N_r$ reactions 
  \begin{equation*}
    \sum_{j=1}^{N_e} \nu_{ij} X_j \leftrightarrows 0, \quad
    i=1,\dotsc, N_r
  \end{equation*}
  where $\nu$ is the stoichiometric matrix. Following
  Morel~\cite{Morel}, we distinguish between \emph{component} and
  \emph{secondary} species by extracting a full rank matrix from
  $\nu$. Component species are a minimal subset of the species such
  that the other secondary species can be written in terms of them (in
  a unique way). Each secondary species gives rise to a reaction that
  expresses how it is formed in terms of the components, and to a mass
  action law that gives the value of its \emph{activity} in terms of
  the component activities. Similarly, each component gives rise to a
  conservation equation, expressing how the given total concentration
  of such a component is distributed among the component itself
  and the secondary species.

Additionally, in the context of reactive transport, it is required to
know how the species are split between those that are in solution, and
those that have been adsorbed on the solid matrix (in this paper we do
not take precipitation into account). We thus introduce (with
obviously $N_e=N_c+N_s+N_x+N_y$)
\begin{itemize}
\item mobile components $c_j, \; j=1,\dots,N_c$, 
\item fixed components $s_j, \;j=1, \dots, N_s$, 
\item mobile secondary species $x_j,\; i=1,\dots, N_x$,
\item fixed secondary species $y_j. \; i=1, \dots, N_y$. 
\end{itemize}
We have identified the name of the species with their concentrations,
and we assume an ideal solution (activities and concentrations are
identified). Mobile secondary species $x$ can be expressed as linear
combinations of mobile components while secondary fixed species depend
on both mobile and fixed components. Therefore the mass action laws
are written as
\begin{equation}
  \label{eq:mal}
  x_i = {K_x}_i \prod_{j=1}^{N_c} c_j^{S_{ij}}, \quad i=1, \dotsc, N_x, \qquad 
  y_i = {K_y}_i \prod_{j=1}^{N_c} c_j^{A_{ij}}\, \prod_{j=1}^{N_s}
  s_j^{B_{ij}}, \quad i=1, \dotsc, N_y,
\end{equation}
where ${K_x}_i$ and ${K_y}_i$ are the equilibrium constants, and
$S_{ij}$, $A_{ij}$ and $B_{ij}$ are the entries of the stoichiometric
matrices $S \in \Rn[N_c \times N_x]$, $A \in \Rn[N_c \times N_y]$ and $B \in
\Rn[N_s \times N_y]$.

Mass conservation for each component is expressed in the form
\begin{equation}
  \label{eq:masconv}
  c + S^T x + A^T y = T, \quad s + B^T y = W,
\end{equation}
where $T_j$ is the total concentration of the mobile component $j$,
and $W_j$ is the total concentration of the fixed component $j$ ($T$
and $W$ are vectors of size $N_c$ and $N_s$ respectively). In the case
of ion exchange, the second mass conservation equation is simply
$B^Ty= W$, and $W$ is the Cationic Exchange Capacity of the porous
matrix (see Appelo and Postma~\cite{AppelPost:05}). As will be seen
later, in the context of coupled transport and chemistry, $T_j$ is
given by the transport model and $W$ is constant. In a closed
chemical system, $T_j$ would be part of the data (total concentration of the
components).

Due to the wildly different orders of magnitude of the concentrations
that are commonly encountered, the chemical problem is reformulated by
using as main unknowns the logarithms of the concentrations. This has
the added advantage that concentrations are automatically positive,
and has become the standard way to solve the
problem~\cite{vanderlee93}. An additional advantage has been pointed
out by Samper et al.~\cite{samxuyang:08}: by taking the logarithms of
the concentrations as unknowns, the Jacobian of the nonlinear system
is symmetric, and with a proper choice of the component species, it
can be shown to be diagonally dominant, and thus nonsingular. The
symmetry can also be seen on equation~\eqref{eq:jacchimie} below.  Let
$\log u$ be the vector with entries $\log u_i$, where $u_i$ are the
entries of vector $u$. Equations~\eqref{eq:mal} can then be rewritten
as a linear system
\begin{equation}  
  \label{eq:mallog}
  \begin{aligned}
    \log x &= \log K_x + S \log c \\
    \log y &= \log K_y + A \log c + B \log s 
  \end{aligned}
\end{equation}

The nonlinear system of equations~\eqref{eq:masconv}
and~\eqref{eq:mallog} forms what will be called the \emph{chemical problem}.
In the sequel, it will be assumed that this problem always has a
(positive) solution $(c, s)$, for all feasible values of the data $T$
and $W$.  This is true in our simplified settings because the chemical
equilibrium problem is a consequence of the minimization of the Gibbs
free energy, which can be shown to be convex in the absence of
minerals (see~\cite{shapshap:65}).

To solve the chemical problem, a variant of Newton's method is used.
As is well known, Newton's method is not always convergent, unless the
initial point is sufficiently close to the solution. However, and this
is especially true in the context of a coupled code where the
chemical problem will be solved repeatedly, it is essential to ensure
that the solver ``never'' fails. We have found that using a
globalized version of Newton's method (using a line search,
cf.~\cite{Kelley95}) was effective in making the algorithm converge
from an arbitrary initial guess. In order to get a smaller system, the
secondary concentrations are eliminated, and the system to be solved involves
only $lc= \log c \in \Rn[N_c]$ and $ls = \log s \in \Rn[N_s]$. Define
the function $H : \Rn[N_c+N_s] \to \Rn[N_c+N_s]$ by
\begin{equation}
\label{eq:chemfun}
  H  \begin{pmatrix} lc \\ ls \end{pmatrix} = 
  \begin{pmatrix}
    \exp(lc) + S^T \exp(\log K_x + S lc) + A^T \exp(\log K_y + A lc +
    B ls) \\
    \exp(ls) + B^T \exp(\log K_y + A lc + B ls), 
  \end{pmatrix}
\end{equation}
where the notation $\exp(v)$ for a vector $v$ means the vector with
elements $\exp(v_j)$, then equations~\eqref{eq:masconv} and
\eqref{eq:mallog} are equivalent to:
\begin{equation}
  \label{eq:chempb}
  H \begin{pmatrix} lc \\ ls   \end{pmatrix} 
  =
  \begin{pmatrix}
    T \\ W
  \end{pmatrix}.
\end{equation}
This is the nonlinear system that to be solved for $lc$ and
$ls$, given $T$ and $W$. The secondary concentrations can then be
computed from equation~\eqref{eq:mallog}.

\medskip
When solving the coupled problem, the distribution of the species
between their mobile form and their fixed form will be needed. The
individual concentrations must still be solved for, but they are
intermediate quantities. Once the component concentrations have been
computed as described in the previous paragraph, one can compute for
each species its mobile part $C_j$ and its fixed part $F_j$ by
\begin{equation}
  \label{eq:defCF}
C=c+S^T x, \qquad F=A^T y.\\
\end{equation}
Note that, by definition, the relationship $T=C+F$ holds. 

In the formulation to be presented below, it will be convenient to
represent the mapping from the vector of total concentrations to the
vector of fixed concentrations. This mapping, denoted by $\Psi$, is
defined by first solving the chemical problem~\eqref{eq:chempb}, then
computing $F$ by~\eqref{eq:defCF}. More precisely
\begin{equation}
  \label{eq:defpsi}
  \begin{aligned}
    \psi : \quad & \mathbf{R}^{Nc} \to \mathbf{R}^{Nc} \\
           & T \mapsto \psi(T) = A^T y,
  \end{aligned}
\end{equation}
where equation~\eqref{eq:chempb} is first solved for $lc$ and $ls$,
then $y$ is computed by~\eqref{eq:mallog}. 

It is important to keep in mind that computing $\Psi(T)$ means solving
the chemical system (plus some simple computations), as this will be the
most expensive part when evaluating the residual of the coupled system
(see eq.\eqref{eq:residual} in section~\ref{sec:NewtKryl}). 

\medskip
As this will be useful later on, the computation of the Jacobian of
$\Psi$ is outlined here. Assume $\Psi(T)$ itself has been computed,
so that the nonlinear system~\eqref{eq:chempb} has been solved. First,
the Jacobian matrix of $H$  should also be computed as
part of the solution process. This is almost certainly needed for
solving the chemical problem, if Newton's method is
used. Differentiating equation~\eqref{eq:chemfun} leads to:
\begin{equation}
  \label{eq:jacchimie}
  H' \begin{pmatrix} lc \\ ls \end{pmatrix} = 
  \begin{pmatrix} \diag(\exp(lc)) & 0 \\ 0 & \diag(\exp(ls))
  \end{pmatrix} +
  \begin{pmatrix} S^T & A^T \\ 0 & B^T \end{pmatrix}
  \begin{pmatrix} \diag(x) & 0\\ 0 & \diag(y) \end{pmatrix}
  \begin{pmatrix} S & 0 \\ A & B \end{pmatrix},
\end{equation}
where $\diag(v)$ is the diagonal matrix with vector $v$ along the diagonal.
Then, by an application of the implicit function theorem (see
for instance~\cite{Rudin:76}), and by differentiating
equation~\eqref{eq:mallog}, there comes
\begin{equation}
\label{eq:psiprim}
  \psi'(T) = A^T \diag(y) \begin{pmatrix}A &B \end{pmatrix}
  \left(H' \begin{pmatrix} lc \\ ls \end{pmatrix}\right)^{-1} 
  \begin{pmatrix} I \\ 0 \end{pmatrix}.
\end{equation}

It should be stressed that the Jacobian of $H$ is needed to computed
the Jacobian of $\Psi$ (inverting it is straightforward, as this will
usually be a small matrix). This may prove problematic in practice for
several reasons. First, the chemical solver may not give access to the
Jacobian, even if it is used internally. This is a limitation to the
``black-box'' approach. Second, for more realistic chemical models,
including non-ideal chemistry, and taking minerals into account,
computing the Jacobian may be much more difficult than the fairly
simple computation outlined above. As a last resort, one could compute
the Jacobian by finite differences, but it will be argued in
section~\ref{sec:NewtKryl} that, for this particular problem, the
analytical computation is more efficient.

\subsection{Coupled transport and chemistry}
\label{sec:couplmodel}

The starting point for the coupled model is the following set of
equations for the total, mobile and fixed concentrations of each
component 
\begin{equation}
\label{sy.1}
\left\{
\begin{aligned}
& \phi\dfrac{\partial C_j}{\partial t} + \phi\dfrac{\partial
  F_j}{\partial t} +L(C_j) = 0 \quad j=1,\dotsc, N_c
\\
& \dfrac{\partial W_j}{\partial t} = 0, \quad j=1,\dotsc, N_s
\end{aligned}
\right.
\end{equation}
These equations can be derived from the individual conservation
equations by standard algebraic manipulations, see for instance Yeh
and Tripathi~\cite{yehtrip}. It is the formulation given in the
benchmark definition~\cite{carray:bench08}, see
also~\cite{saalayorcarr98},~\cite{dieulerhkern:09}. The second
equation is obvious, as $W_j$ was taken as a constant (at each point
in space).

Taking into account the relation $T_j=C_j+F_j,\; j=1,\dotsc, N_c$
noted above, the first equation of the system is equivalent to 
\begin{equation}
\label{eq:transcoupl}
\phi\dfrac{\partial T_j}{\partial t} +L(C_j) = 0\quad j=1,\dots,N_c, 
\end{equation}
where $T_j$ is the total concentration, $C_j$
the total mobile concentration, and $F_j$ the total
fixed concentration for component $j$.

From now on, $L$ will denote the discretized transport operator, as
defined in equation~\eqref{eq:trans_ode}. Each
unknown concentration depends on both the grid point index, and the
chemical species index. We will use a notation inspired from Matlab.
For a concentration $u_{ij}$, where $i \in [1, N_g]$ represents the
spatial index and $j \in [1, N_c]$ represents the chemical index, we
shall denote by
\begin{itemize}
\item $u_{:,j}$ the column vector of concentrations of species $j$ at all
  grid points;
\item $u_{i,:}$ the row vector of concentrations of all chemical species
  in grid cell $x_i$. 
\end{itemize}
The unknowns will be numbered first by chemical species, then by
grid points. Thus all the unknowns for a single grid point are numbered
contiguously. 

The coupled problem is obtained by putting together
equation~\eqref{eq:transcoupl} above with the definition of the
chemical solution operator $\psi$, defined in eq.~\eqref{eq:defpsi}
(the subscript T denotes transposition):
\begin{equation}
  \label{eq:syscoupl}
  \left\{
    \begin{aligned}
      & M\dfrac{\partial C_{:,j}}{\partial t}+ M\dfrac{\partial
        F_{:,j}}{\partial t} +L(C_{:,j}) = g_{:,j} & &
      \; j=1,\dots,N_c \\
      & T_{ij} = C_{ij}+F_{ij}, & i=1, \dots, N_g, & \; j=1, \dots, N_c \\
      & F_{i,:} = \psi(T_{i,:}^T)^T, & i=1, \dots, N_g &
    \end{aligned}
  \right.
\end{equation}
This system is then discretized in time to obtain the fully discrete
coupled nonlinear system. In this work we restrict to a simple
backward Euler scheme with constant step--size, noting that other
more sophisticated, strategies are obviously possible (in particular,
an adaptive step-size is essential for efficiency). Denoting time
indexes by a superscript, the following system is obtained
\begin{equation}
  \label{eq:sysnl}  
  \left\{
    \begin{aligned}
      & M \dfrac{C^{n+1}_{:,j} - C^n_{:,j}}{\Delta t} + 
      M \dfrac{F^{n+1}_{:,j} - F^n_{:,j}}{\Delta t} +
      L(C^{n+1}_{:,j} ) = g_{:,j}  & & \; j=1,\dots,N_c \\
      & T^{n+1}_{ij}=C^{n+1}_{ij}+F^{n+1}_{ij} &  i=1, \dots, N_g, &\;
      j=1, \dots, N_c \\
      & F^{n+1}_{i,:}=\psi((T^{n+1}_{i,:})^T)^T & i=1, \dots, N_g &
    \end{aligned}
  \right.
\end{equation}

This is the system to be solved at each time step.  

\section{Formulation and coupling algorithms}
\label{sec:couplalg}
The formulation of reactive transport seen above gives rise to a large
system of nonlinear equations. For complex problems, its solution will
require a large amount of computer time, which makes it important to
choose an appropriate method. In this section, several formulations
and approaches that have appeared in the literature will be reviewed.

Thanks to the relationship $T=C+F$, it is easy to eliminate one of the
3 variables, and this leads to different formulations for the coupled
problem, depending on which variables are kept in the transport
equation. We keep the system continuous in time, as it makes the
notation somewhat lighter, but the same manipulations can obviously be
done at the discrete level too.

According to Saaltink et al.~\cite{saalayorcarr98}, see also
Salignac~\cite{sali:98}, one can derive three main formulations from
the system given in~\eqref{eq:syscoupl}:
\begin{itemize}
\item formulation (TC) where $T$ is the principal variable , $C$ the
  transported variable
\begin{equation}\label{co.3} M\dfrac{\partial T_{:j}}{\partial
    t} + L(C_{:j})=g_{:,j}\end{equation}
This is the formulation used by Erhel et
al. in~\cite{dieulerh07,dieulerhkern:09}, as it 
lends itself best to a DAE type algorithm. It is not convenient for
our purpose, as the transport equation then involves both $T$ and $C$,
and is thus not easily used with an existing transport solver.
\item  formulation (TT) where $T$ is the principal variable , $T$
  transported variable 
\begin{equation}\label{co.4} 
  M\dfrac{\partial T_{:j}}{\partial t} + L(T_{:j})+L(F_{:j})=g_{:,j} 
\end{equation} 
This seems to be the least satisfactory formulation, as the transport
operates on the fixed species, and for this reason it will not be considered further. 
\item  formulation (CC) where $C$ is the principal variable , $C$
  transported variable 
\begin{equation}\label{co.5} 
M\dfrac{\partial C_{:j}}{\partial
    t}+M\dfrac{\partial F_{:j}}{\partial t} + L(C_{:j})=g_{:,j} 
\end{equation}
This is formulation 4 in Saaltink et al.~\cite{saalayorcarr98}, and is
the formulation chosen below. It has been reported that this
formulation is the least suitable for use in an operator split
algorithm, because $C$ and $F$ are used at different time levels (to
compute the data for the chemical problem). When this formulation is
used in a global method this should not matter as much, as the
iterations are ran to convergence, and both values should eventually
get close to their limits.
\end{itemize}

Formulation (CC) will be used in the rest of the paper,
because it takes the form of a standard transport operator, with a
source term coming from the chemical part. Its structure is closely
related to the system describing single species transport with
sorption, as seen for instance in~\cite{KacFrol:06}, or~\cite{kamike:03},
with the main differences that the unknown is a vector
of concentration, and mostly that what plays the role of the sorption
isotherm is the implicitly defined function $\Psi$ introduced
in~\eqref{eq:defpsi}.

\subsection{Review of former approaches}

At each time step, the system given by~\eqref{eq:sysnl} (one transport
equation for each component and one chemical system for each grid
point) forms a large nonlinear system, whose size is the number of
components times the number of grid points. This system has
traditionally been solved by a sequential two-steps approach, as
reviewed below (cf.~\cite{yehtrip}).  However, this method suffers
from several defects: it may severely restrict the step size to ensure
convergence, and if used non-iteratively it is only first order in
time, which may introduce additional errors
(cf.~\cite{carrmosbeh:04}). Due to its quadratic convergence rate,
Newton's method would be an ideal candidate for solving the system. On
the other hand, a practical difficulty has to be reckoned with:
Newton's method requires the solution of a linear system with the
Jacobian matrix at each iteration step. In realistic situations, it
will not be possible to store, much less factor, the Jacobian matrix.
As will be seen in section~\ref{sec:NewtKryl}, this difficulty can be
overcome by resorting to an iterative method for solving the linear
system.

\subsubsection{Sequential approach}
\label{sec:seqcoupl}

The sequential approach consists of separately solving the chemical
equations and the transport equations. The method has been used in
numerous papers: see for instance~\cite{yehtrip}, and
also~\cite{kamike:03},~\cite{LuBuOl:00},~\cite{saalayorcarr98},~\cite{carrmosbeh:04}
or~\cite{maymacq:09}<.
At each iteration, a transport equation for each component is solved
first, with a source term given by the (change in) fixed concentration
at the previous iteration. This total mobile concentrations will be
added to a total fixed concentration computed in the previous
iteration, to obtain the total used as data for solving a chemical
problem at each grid point. These steps are then iterated until
convergence.

In the geochemical literature, this is known as an operator splitting
approach (usually called Standard Iterative Approach, or SIA), but it
is more properly a block Gauss-Seidel methods on the coupled system,
as each subsystem is solved alternatively. The method is quite
appealing, as it is easy to implement starting from separate transport
and chemistry codes, and can provide good accuracy if implemented
carefully, as shown in the references above). As will be seen below,
theses advantages can be retained in the Newton--Krylov framework.

The Standard Non-Iterative Approach (SNIA) is the case where only one
iteration of the method is carried out at each time step. In that
case, splitting errors can become important, and the method is not
really suitable for difficult problems. 

The SIA approach does not suffer from splitting errors if the
tolerance is small enough, but it may require a small time step to
obtain convergence in the case of stiff problems. The main drawback of
the method is thus that the size of the time step is used to control
convergence, and not based on the physical character of the solution.

\subsubsection{Direct Substitution Approach}
\label{sec:DSA}

As computing power increased, it was recognized that the operator
splitting methods of the previous sections could not satisfactorily
handle difficult problems, and more tightly coupled method came to
more widespread use. 

The Direct Substitution Approach method consists in solving for the
individual concentrations of the components, that is
\emph{substituting} equations~\eqref{eq:masconv}--\eqref{eq:mallog} in
equation~(\ref{eq:transp}) (this can be done explicitly, as in Hammond
et al.~\cite{hammvallicht05}, or implicitly, as in Kräutle et
al.~\cite{knabner05,knabner07}, or Saaltink et
al.~\cite{saalayorcarr98}).  It is also possible to reformulate the
problem as a differential algebraic system (DAE), and to take
advantage of the high quality software available for such problems, as
in Erhel et al.~\cite{dieulerhkern:09},~\cite{dedieuerh:09}
or~\cite{dieulthes:08} . A high performance parallel implementation is
described by Hammond et al.~\cite{hammvallicht05}, using a
Jacobian--Free Newton--Krylov method (see section~\ref{sec:NewtKryl}).

The main advantages of this approach are to avoid the errors caused by
the separation of operators, and to allow fast convergence
independently of the time step, but its principal drawback is the need
to form and to store the Jacobian matrix especially for a large
problem. Moreover, sometimes it may be difficult to calculate the
exact derivatives for geochemical processes especially when
precipitation phenomena or kinetic reactions are taken into account.

The size of the system can be made smaller by means of a reduction
method, cf Kräutle et al.~\cite{knabner05,knabner07},
and~\cite{HKP:09}. The reduction method makes a
change of variables in the chemical system, so that a set of decoupled
transport equation is first solved, leaving a smaller nonlinear
system, that is still solved with Newton's method.

\subsection{A Newton--Krylov based fully coupled method}
\label{sec:NewtKryl}
As was already mentioned in the previous section, Hammond et
al.~\cite{hammvallicht05} have used a Newton--Krylov method for
solving the system obtained from the DSA approach. Substituting the
chemical equations in the transport operator is the most
straightforward way of formulating the coupled problem, but leads to a
system where chemistry and transport terms are mixed, and makes it
virtually impossible to separate the transport and chemistry
modules. However, this separation is seen as one of the important
advantages of the operator splitting approaches. 

By coupling the formulation given in section~\ref{sec:couplmodel} with
the Newton--Krylov framework, a strongly coupled method that can be
implemented by keeping transport and chemistry separate is obtained.
Thus, the chemical equations are not directly substituted in the
transport equation, but the function $\Psi$ introduced previously
in~\eqref{eq:defpsi} is used to represent the effect of chemistry.
Different formulations could be adopted depending on the choice of
unknowns (refer back to section~\ref{sec:couplalg}). In this work,
both the total mobile and fixed concentrations, and also the total
concentrations (though they could easily be eliminated) are chosen as
main unknowns. 

Even though this method may be more expensive than the methods based
on DSA, its main advantage is to make it possible to treat chemistry
as a black--box, even in the Newton--Krylov context. This may be
important, as chemical simulators are becoming increasingly
sophisticated.

Recall (equation~\eqref{eq:sysnl}) that the nonlinear
system to be solved at each time step is
\begin{equation}
  \label{eq:sysnl1}  
  \left\{
    \begin{aligned}
      & (M+\Delta t L) C^{n+1}_{:,j} + M F^{n+1}_{:,j} -
      b_{:,j}^n & = 0, & \quad j=1, \dots, N_c \\ 
      & T^{n+1} - C^{n+1} - F^{n+1} &= 0,  & \\
      & F^{n+1}_{i,:} - \psi\left((T^{n+1}_{i,:})^T\right)^T &= 0,  &
      \quad i=1,  \dots, N_g 
    \end{aligned}
  \right.
\end{equation}
where $b_{:,j}^n =M C^{n}_{:,j}+\Delta t\;g_{:,j}^{n+1} + M \Delta
 t F^{n}_{:,j}$ is known.

Denoting by $G: \Rn[3 N_c N_g] \to \Rn[3 N_C N_g]$ the function
\begin{equation}
\label{eq:residual}
  G  \begin{pmatrix}  C \\ T \\ F \end{pmatrix} = 
  \begin{pmatrix}
    \left((M+\Delta t L) C_{:,j} + M F_{:,j} -
      b_{:,j}^n \right)_{j \in [1, N_c]} \\
     T - C - F \\
     \left(F_{i,:} - \psi\left((T_{i,:})^T\right)^T \right)_{i \in [1, N_g]}
  \end{pmatrix},
\end{equation}
the nonlinear problem to be solved at each time step is $G(Z)=0$,
where $Z$ denotes the vector $\begin{pmatrix} C^{n+1} \\ T^{n+1} \\
  F^{n+1} \end{pmatrix}$.

Recall that at each step of the ``pure'' form of Newton's method
for solving $G(Z)=0$, one should compute the Jacobian matrix $J
=G'(Z^k)$, solve the linear system (usually by Gaussian elimination)
\begin{equation}
  \label{eq:newtres}
  J\, \delta Z = - G(Z^k)
\end{equation}
and then set $Z^{k+1} = Z^k + \delta Z$. In practice, one should use
some form of globalization procedure in order to ensure convergence
from an arbitrary starting point. If a line search is used, the last
step should be replaced by $Z^{k+1} = \delta Z + \lambda Z^k$, where
$\lambda$ is determined by the line search procedure. 

The main drawback of the method for large scale problems is again the
need to form, and then factor, the Jacobian matrix. For coupled
problem such as the one studied in this paper, there is the additional
difficulty of simply computing the Jacobian: the numerical methods for
transport and chemistry are quite different, and it is even possible
that the simulation codes have been written by different groups.

The Newton--Krylov method (see~\cite{Kelley95},~\cite{knollkeyes04}
and~\cite{hammvallicht05}, to which our work is closely related), is a
variant of Newton's method 
where the linear system that arises at each step of Newton's method is
solved by an \emph{iterative} method (of Krylov type). The main
advantage of this type of method is that the full Jacobian is not
needed, one just needs to be able to compute the product of the
Jacobian with a vector. As this is a directional derivative, this
leads to the Jacobian free methods, where this product is approximated
by finite differences. However, for some problems, it may be possible
to compute the needed directional derivative exactly. As will be seen
below, this is the case for our coupled problem, provided 
the Jacobian of the chemical problem can be computed. This is both cheaper
and more accurate. 

The main contribution of this paper is to show that the formulation
given above lends itself to an implementation of Newton's method that
allows to keep the two codes separate. This is in keeping with thte
philosophy set forth in the review paper by Keyes and
Knoll~\cite{knollkeyes04} that a Newton-Krylov solver can often be
made by wrapping a classical split-step solver. This is what is being
done here, as the formulation to which the Newton-Krylov method is
applied is the one used for operator splitting. Additionally, it will
be shown below that the Jacobian may even be formed in block form,
provided the individual codes provide their Jacobians (this is
obviously easier for transport than for chemistry), and this obviously
carries over to the Jacobian--vector product.

At this point, it is appropriate to add a few comments on the size of
the problems envisioned. The examples used in this work are small
scale, one dimensional, problems. They can hardly be called large. On
the other hand, we believe they are representative of the problems
that will be encountered in more realistic applications. For such
problems, in 2 or 3 space dimensions, involving tens or hundreds of
thousands of grid points and several tens of chemical species, the
nonlinear system will indeed be very large, and a method like that of
Hammond et al.~\cite{hammvallicht05}, or like the method presented in
this section will be necessary.

\medskip
A Krylov subspace method (see for instance~\cite{Kelley95}) is used to
approximately solve the linear system in equation~(\ref{eq:newtres}). The
linear iterates are drawn from the Krylov subspace, $K_j = \text{span}
\{r_0,Jr_0,J^{2}r_0, \dotsc ,J^{j-1}r_0\}$. In the GMRES method
(see~\cite{sasc:86}), the iterates are defined to minimize the residual
${\|J\delta Z_{j}+G(Z)\|}_2$ over $K_j$. Other methods, such as
Bi-CGSTAB~\cite{vanderVorst:1992:BCF} or QMR~\cite{frna:91a} could be
used as well. 

As the linear system is not solved exactly, the convergence theory for
Newron's method does not apply directly. However, the theory has been
extended by Dembo et al.~\cite{dembo:400} to the class of Inexact Newton methods,
of which the Newton--Krylov methods are representatives. The main
consequences of this analysis are summarized below.

An important issue in such methods is the stopping criterion for the
inner linear iteration. A stopping criterion of the form
\begin{equation}
  \| J \delta Z + G(Z^k) \| \le \eta_k \| G(Z^k) \|
\end{equation}
in this context, as the initial iterate is usually 0. The choice of
the \emph{forcing term} $\eta_k$ should strike a balance between two
conflicting goals:
\begin{itemize}
\item Keep the (local) convergence of Newton's methods;
\item Avoid over-solving, that is taking too many linear iterations
  when still far away from the nonlinear solution.
\end{itemize}
The first goal will tend to require a small value for $\eta_k$, while
the second one obviously tends to make $\eta_k$ larger. It has been
shown (see theorem 6.1.4 in~\cite{Kelley95}) that provided $\eta_k$ is
bounded away from 1, the inexact Newton's method will converge, and
that superlinear convergence obtains if $\eta_k$ goes to zero faster
than $\| G(Z^k)\|$.  Based on this result, the strategy proposed by
Kelley in~\cite{Kelley95} (after the choice in~\cite{eisenwalk:96})
computes $\eta_k$ as
\begin{equation}
\label{eq:forcing}
  \eta_k = \gamma\, {\|G(Z_k) \|^2} / {\| G(Z_{k-1})\|^2},
\end{equation}
where $\gamma \in [0, 1]$ is a parameter (the value suggested
in~\cite{Kelley95} is $\gamma=0.9$). Safeguards are added to this
choice in order to prevent $\eta_k$ to become too close to 1, or too
small. It is also necessary to globalize the algorithm, and this can
be done using a line search, just as in the ``classical'' Newton's
method.

The other main practical advantage of the Newton--Krylov methods is
that they do not require forming the Jacobian matrix. All that is
needed is the ability to compute the product of the Jacobian matrix by
an arbitrary vector, in order to enlarge the Krylov subspace. This
matrix--vector product can be interpreted as a directional derivative.
This means that, for complex functions $G$ it may not be necessary to
compute the Jacobian, at the cost of one extra evaluation of the
function itself. It turns out, however, that in our case, this
trade-off is not advantageous. Indeed, it is well known that the most
expensive part of the evaluation of $G$ is the solution of the
chemical problem at each grid point. On the other hand, it was shown
above that computing the Jacobian of $\psi$ is actually cheaper than
computing $\psi$ itself (once $\psi$ has already been computed), as it
only involves the solution of a linear system (see
equation~\eqref{eq:psiprim}), whereas computing $\psi$ itself requires
the solution of a nonlinear system.

\medskip
It will now be shown how the method can be implemented, given modules
for transport and chemistry. 

The first ingredient needed is the computation of the residual, that
is evaluating the function $G$ defined in~\eqref{eq:residual}. Given a
vector $Z= \left( \begin{smallmatrix}  C\\ T \\ F \end{smallmatrix} \right)$, $Z$ is
first split into its three components, and each sub-vector is regarded
as a $N_g \times N_c$ matrix, as in section~\ref{sec:couplmodel}. Then
$G(Z)$ is computed by block:
\begin{itemize}
\item For the transport block, the transport operator is applied to
  each species $C_{:,j}$, with a source term given by $-M \dfrac{F_{:,j}
    - F^n_{:,j}}{\Delta t}$, for $j=1, \dotsc, N_c$
  ($F^n$ denotes $F$ at the previous time step);
\item The second block is the trivial computation $T-C-F$;
\item The third block is the solution of the chemical problem at each
  grid point: $F_{:,i}-\Psi(T_{i,:})$, for $i=1, \dotsc, Ng$. 
\end{itemize}
This shows that the first block will only need transport related
quantities, whereas the third block will only call chemistry related
ones. Actually, these are the same computations that would be needed
for implementing a operator splitting method. 

As far as the Jacobian matrix--vector product is concerned, and using
the computation in section~\ref{sec:chimie}, the action of the
Jacobian on a vector $ v = \left(
\begin{smallmatrix}  v_c\\ v_T \\v_F \end{smallmatrix} \right)$ (that
is the directional derivative of $G$ in the direction of the vector $v$)
can be computed as
\begin{equation}
\label{eq:jacvec}
  J  \begin{pmatrix} v_C \\ v_T \\ v_F \\ \end{pmatrix} = 
  \begin{pmatrix}
    \bigl( (M + \Delta t L){v_C}_{:,j} + M {v_F}_{:,j}\bigr)_{j
      \in [1, N_c]} \\[1.2ex] 
    -v_C + v_T - v_F \\[1.2ex]
    \bigl( {v_F}_{i,:} - {v_T}_{i,:} (\psi'(T_{i,:}^T))^T  \bigr)_{i
      \in [1, n_x]} 
  \end{pmatrix}.
\end{equation}

Even though it is not used as such in this work, it is valuable to
examine the structure of the Jacobian. As the previous computation
shows, the Jacobian also has a natural block structure. Recall that
the unknowns are numbered by species at each point in space. Then the
block corresponding to the action of $L$ can be written using the
Kronecker product (see for instance~\cite{Horn:90}) as $A= (M+\Delta t
L) \otimes I$. Then the Jacobian matrix is
\begin{equation}
  \label{eq:jacstored}
J=
\begin{pmatrix}
A & 0& M \\
 -I & I & -I\\
 0 & -\psi'(T^T) & I 
\end{pmatrix},
  \end{equation}
where $\psi'(T) = \text{diag}(\psi'(T_{1,:}^T), \dots, \psi'(T_{N_g,:}^T)$
is the Jacobian of $\psi$, and for each $i=1, \dotsc, N_g$,
$\psi'(T_{i,:}^T)$ is a small $N_c$ by $N_c$ 
block.  
\begin{figure}[thb]
  \begin{center}
    \includegraphics[width=0.5\textwidth,draft=false]{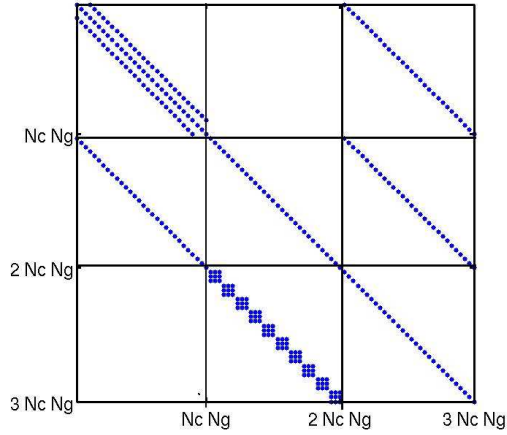}
  \end{center}
  \caption{The block structure of the Jacobian matrix}
  \label{fig:jac}
\end{figure}
The structure of the Jacobian is illustrated on figure~\ref{fig:jac},
for the case $N_g=10$, $N_c=3$. It is a $3\times 3$ block matrix,
each bock being of size $N_g \times N_c$. We can clearly see the different
parts of the Jacobian: the transport part in the upper left corner has
3 diagonals corresponding to the Kronecker product structure (remember
that $L$ is tridiagonal), and the chemistry part at the bottom has 10
$3 \times 3$ blocks.

It would in principle possible to compute and store the Jacobian matrix
according to equation~(\ref{eq:jacstored}) as a sparse matrix, and to compute the
matrix--vector product using a general purpose routine. The advantage
of the method given in equation~(\ref{eq:jacvec}) is that the
structure of the Jacobian is fully exploited, which leads to a much
more economical computation.

\section{Numerical results}
\label{sec:numres}

\subsection{Ion exchange}
\label{sec:test}
The following example of advective transport in the presence of cation
exchanger is adopted as a first test case comparison of both
approaches. The example is used in the documentation of
PHREEQC-2~\cite{Appelo} as Example~11.  The one-dimensional simulation
problem describes a column experiment where the chemical composition
of the effluent from a column containing a cation exchanger is
simulated. Initially, the column contains a sodium-potassium-nitrate
solution in equilibrium with the cation exchanger. The column is then
flushed with three pore volumes of calcium chloride solution, so that
an equilibrium state with calcium and chloride is reached. Calcium,
potassium, and sodium react to equilibrium with the exchanger at all
times. The flow and transport parameters used for this example are
presented in Table~\ref{tab:2}, and the initial and injected
concentrations are listed in Table~\ref{table1}. The Cationic Exchange
Capacity for the exchanger is $1.1\, \text{mmol}/\text{l}$.

\begin{minipage}[c]{0.6\textwidth}
  \begin{table}[H]
    \begin{tabular}{|l|l|} 
\hline
Darcy velocity & $ 2.78\,10^{-6}\; \text{m}/\text{s}$ \\
Diffusion coefficient & $5.56\, 10^{-9}\; \text{m}^2/ \text{s}$  \\
Length of column & $0.08 \; \text{m}$ \\
Mesh size & $0.0002\; \text{m}$\\
Duration of experiment & 1 day \\
Time step & $720\; \text{s}$ \\
\hline
    \end{tabular}
    \caption{Physical parameters}
    \label{tab:2}
  \end{table}
\end{minipage}
%%\hspace*{-2cm}
\begin{minipage}[c]{0.4\textwidth}
\begin{table}[H]
\begin{tabular}{|l|l|l|}
\hline
Component \T \B & $C_{\text{init}}$ & $C_\text{{inflow}}$ \\
\hline
Ca \T & 0 & $0.6\,10^{-3}$ \\
Cl & 0 & $1.2\,10^{-3}$\\
K \B & $2.0\, 10^{-4}$ & 0 \\
Na \B & $1.0\,10^{-3}$ & 0\\
\hline
\end{tabular}
\caption{Initial and injected concentrations}
\label{table1}
\end{table}
\end{minipage}
The chemical reactions for this example are:
%\vspace*{-1.5cm}
\begin{eqnarray*}
\text{Na}^+ + \text{X}^- & \rightleftharpoons & \text{NaX} \\ 
\text{K}^+ + \text{X}^- &\rightleftharpoons & \text{KX} \\
\dfrac{1}{2}\text{Ca}^{2+} + \text{X}^- & \rightleftharpoons & \dfrac{1}{2}\text{CaX}_2 
\end{eqnarray*}
with $\text{NaX}$, $\text{KX}$ and $\text{CaX}_2$ are (sorbed)
complexes, and $X$ indicates exchange site with charge -1

\subsubsection{Comparison with Phreeqc}

\begin{figure}[thb]
  \centering
  \includegraphics[width=0.45\textwidth]{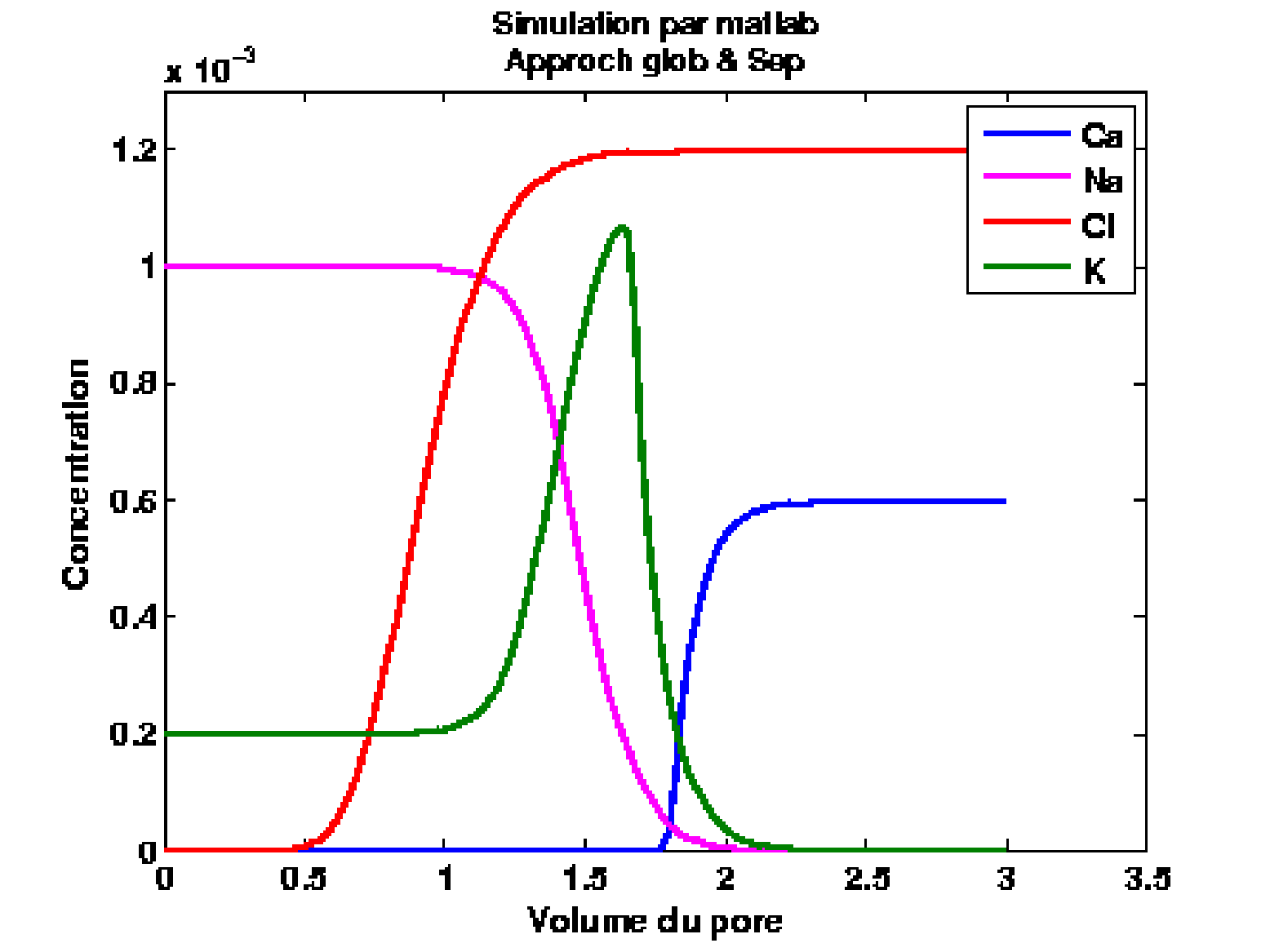}
  \includegraphics[width=0.45\textwidth]{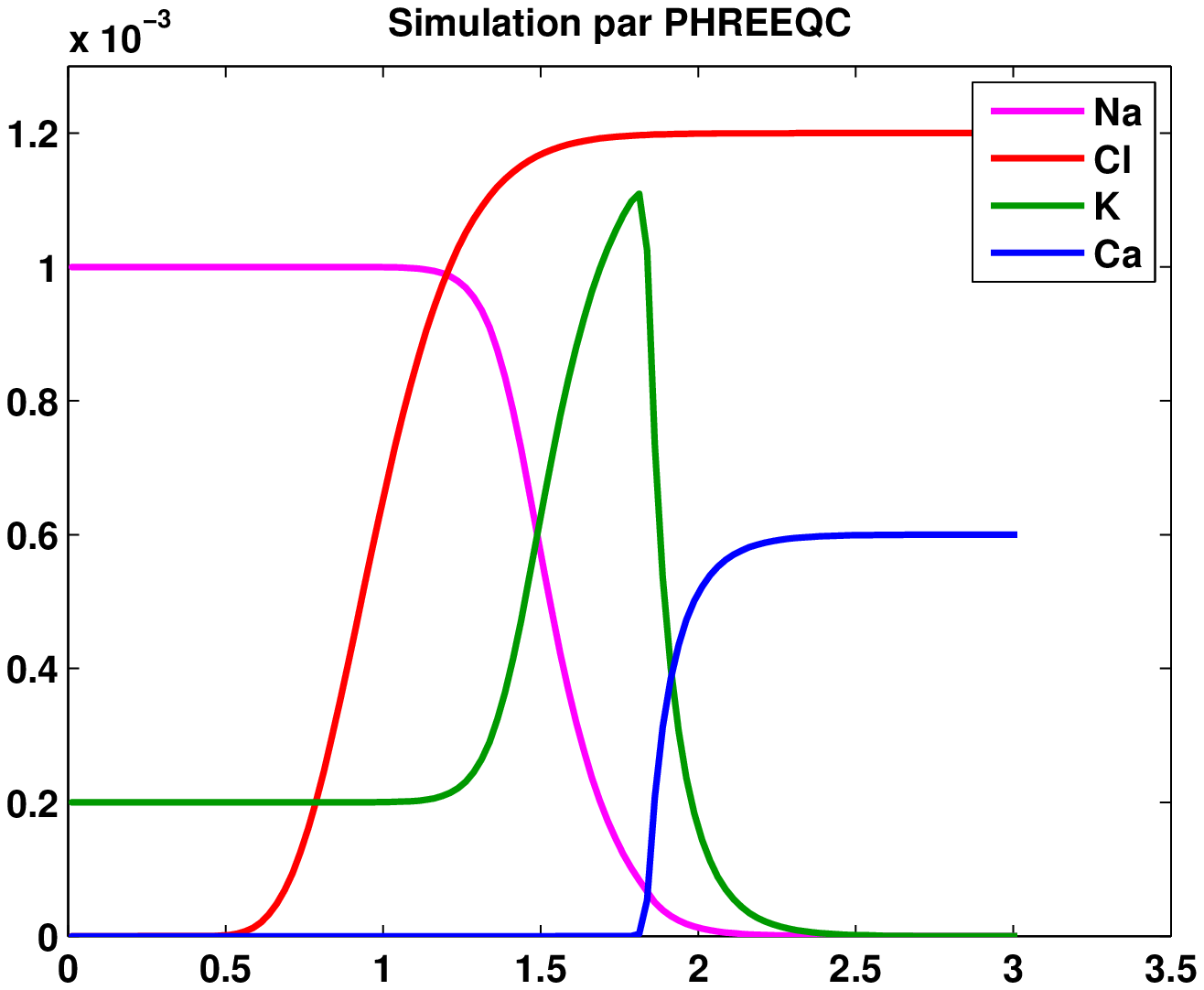}
  \caption{Elution curves (concentrations at the end of the column)
    versus time, for the problem of section~\ref{sec:test}. Left:
    global method, right: PhreeqC reference.} 
  \label{fig:concex11}
\end{figure}

Figure~\ref{fig:concex11} shows elution curves, that is the evolution of
the concentration of the various species at the end of the column, as
a function of time.  The sorbed potassium and sodium ions are
successively replaced by calcium. Because potassium exchanges more
strongly than sodium (as indicated by a larger value of log K in the
exchange reaction), sodium is released first, followed by potassium.
Finally when all of concentration has been released, the concentration
of calcium increases to its steady-state value, the potassium is
displaced from the exchanger and the concentration in solution
increases to balance the $\text{Cl}^{-}$ concentration.

Both the sequential method and the global method described in
section~\ref{sec:NewtKryl} have been applied to the test case
described in section~\ref{sec:test}. Both the computational demands
and the accuracy of the solutions will be compared.  

As can be seen on figure~\ref{fig:concex11}, the results obtained are
close to those computed by PhreeqC. One can still see differences both
in the location and amplitude of the peak in potassium concnetration,
and in the region where the three curves cross. These results are
also comparable to those obtained by Xu et
al.~\cite{XuSamAyoManCus:99}.

\subsubsection{Performance of the method}

\begin{figure}[ht]
  \centering
  \includegraphics[width=0.4\textwidth]{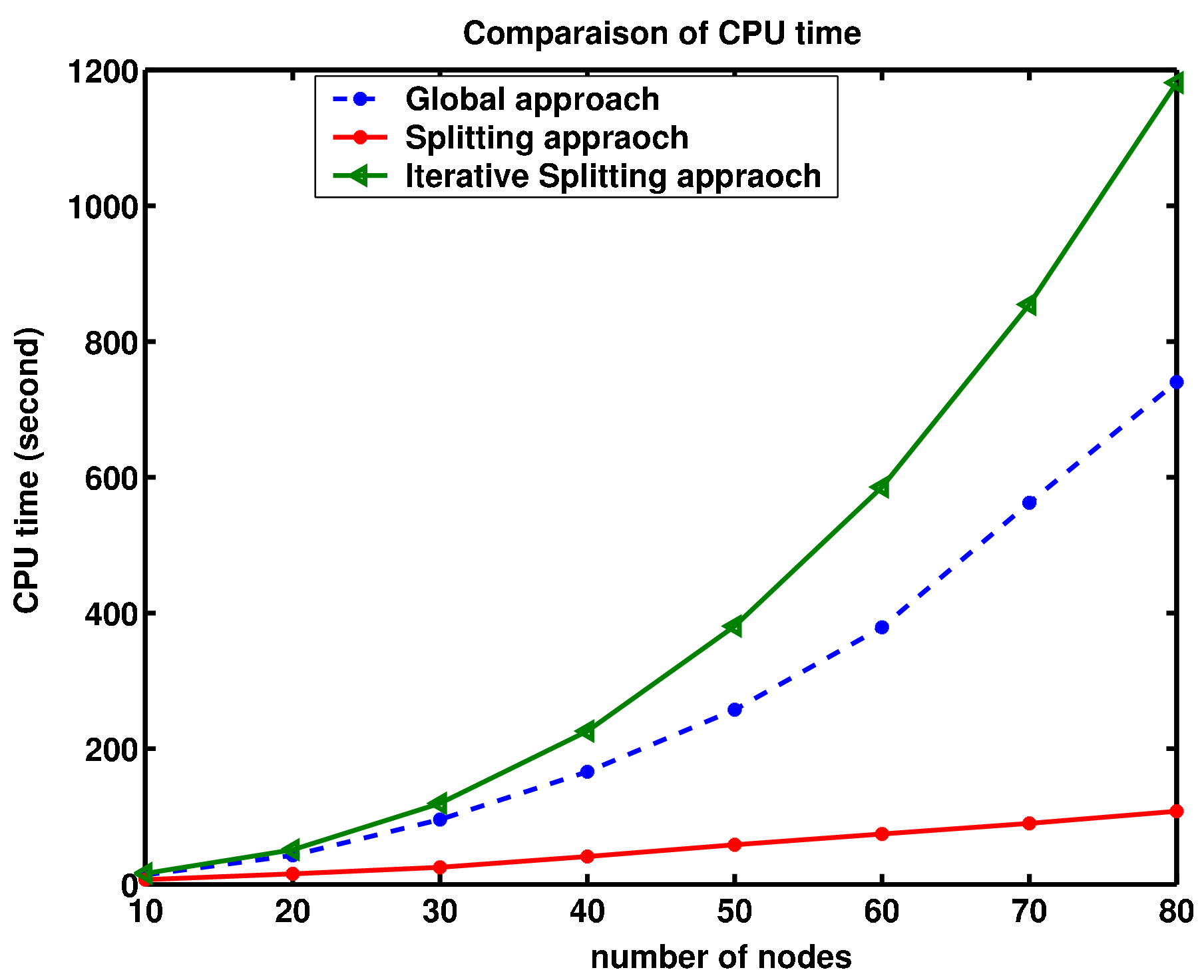}  
  \caption{Computing time for 3 methods applied to the ion exchange of section~\ref{sec:test}}
  \label{fig:cpuex11}
\end{figure}

The CPU times for the iterative splitting, non iterative splitting and
global approaches are compared on figure~\ref{fig:cpuex11}. The CPU
time required for each method is plotted versus the number of the
nodes of the grid. As expected, it can be seen that the
non-iterative method requires  much less CPU time than the iterative
methods. On the other hand, the global approach described in the paper
requires less time than the iterative splitting, at least for the
simple chemical system considered here.  

For a single time step, the iterative splitting approach requires
between 20 and 27 iterations on the average.  The number of fixed
point iterations increase with the number of the nodes in the grid.
On the other hand, the number of Newton iterations for the
Newton--Krylov method is less than 6, independently of the number of
nodes. The number of Krylov iteration for each Newton step, however,
does increase with the number of nodes. We go back to this issue in
subsection~\ref{sec:perfs}.

\subsection{The 1D ``easy'' MoMaS Benchmark}
The global and the splitting approaches will now be applied to the 1D
easy GDR Momas Benchmark, as described in the introductory paper to
this special issue~\cite{carray:bench08}, see also the original
description in~\cite{carray:benchweb}. Let us just recall that the
model is a one-dimensional column, made of 2 different media: the part
in the middle is less conductive but more reactive than the
surrounding medium. The chemical system has 5 components (4 mobile
components and a fixed component), and 7 secondary species. The
equilibrium constants vary over 50 orders of magnitude, and the
stoichiometric coefficients can be as large as 4, making the problem
highly non-linear. 

\begin{figure}[thb]
  \centering
  \subfloat[{$t=10$}]{\includegraphics[width=0.48\textwidth]{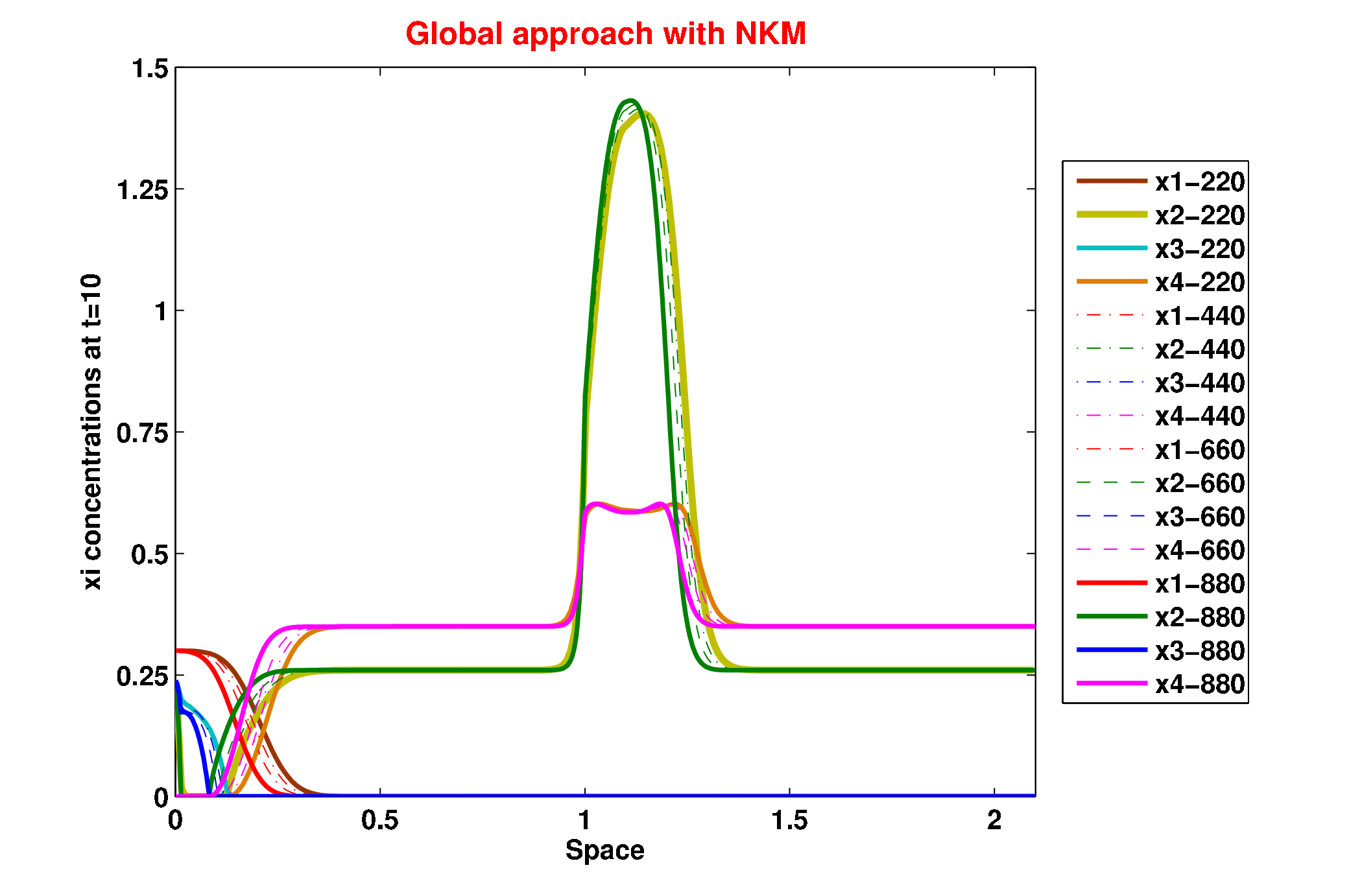}\label{fig:x1x4L}}   
  \subfloat[{$t=50$}]{\includegraphics[width=0.48\textwidth]{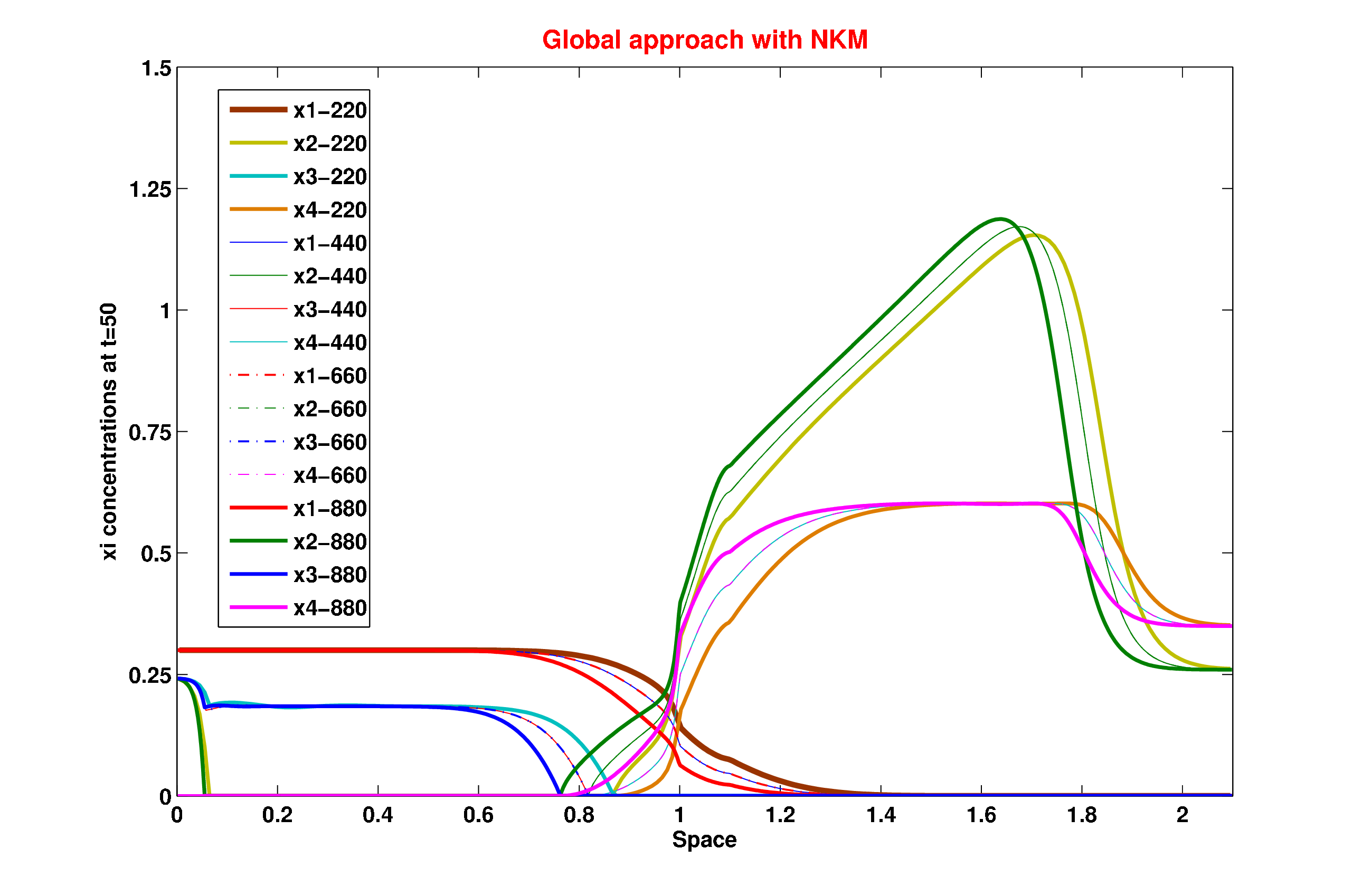}\label{fig:x1x4R}}    
  \caption{Concentration of all components at times $t=10$ and $t=50$,
  for various mesh resolutions}
  \label{fig:x1x4}
\end{figure}

\begin{figure}[h]
  \centering
  \subfloat[$t=0$ to $t=400$]{
    \includegraphics[width=0.48\textwidth]{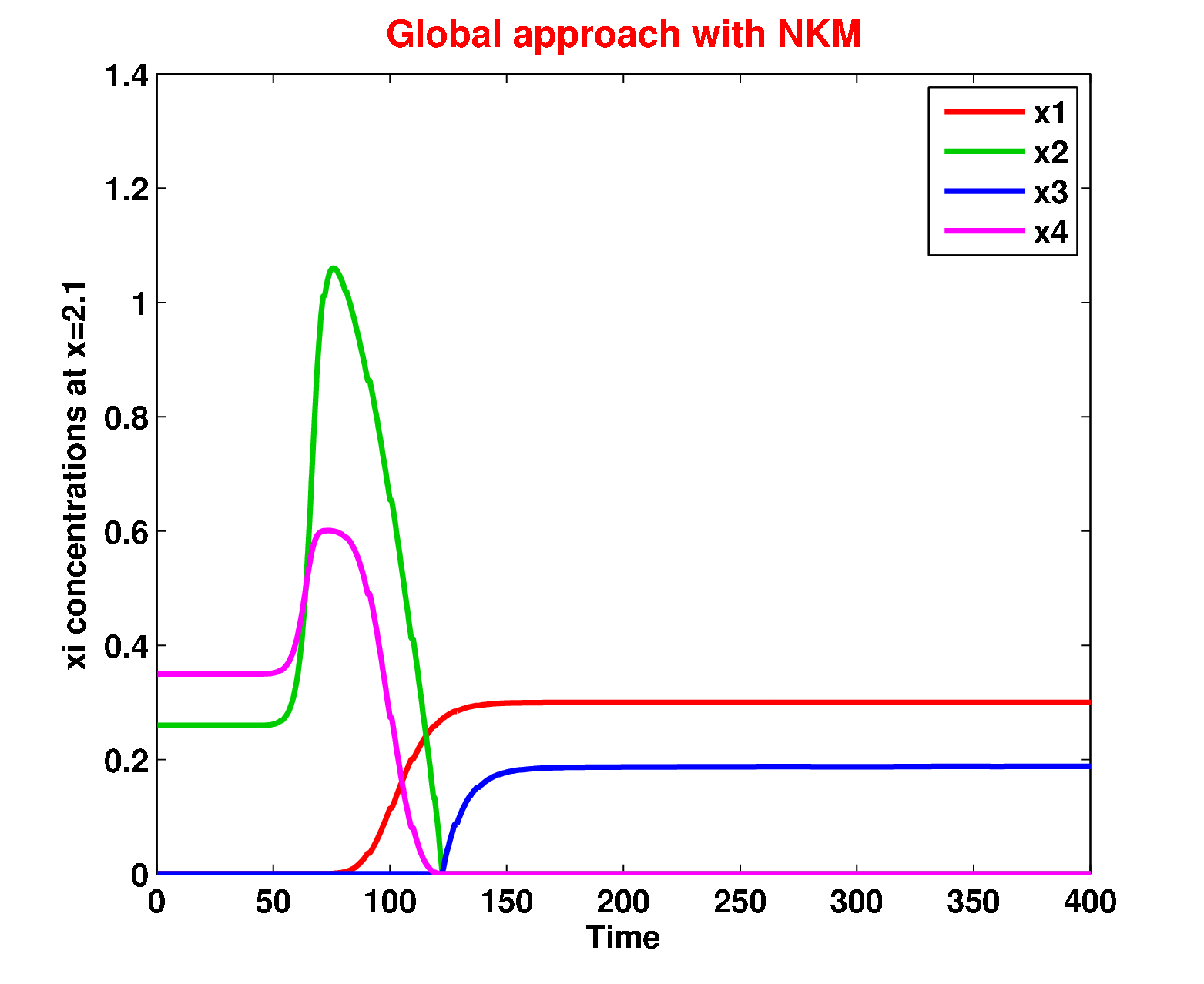}\label{fig:x1x4tpsL}
}
\subfloat[$t=4900$ to $t=5300$]{
  \includegraphics[width=0.48\textwidth]{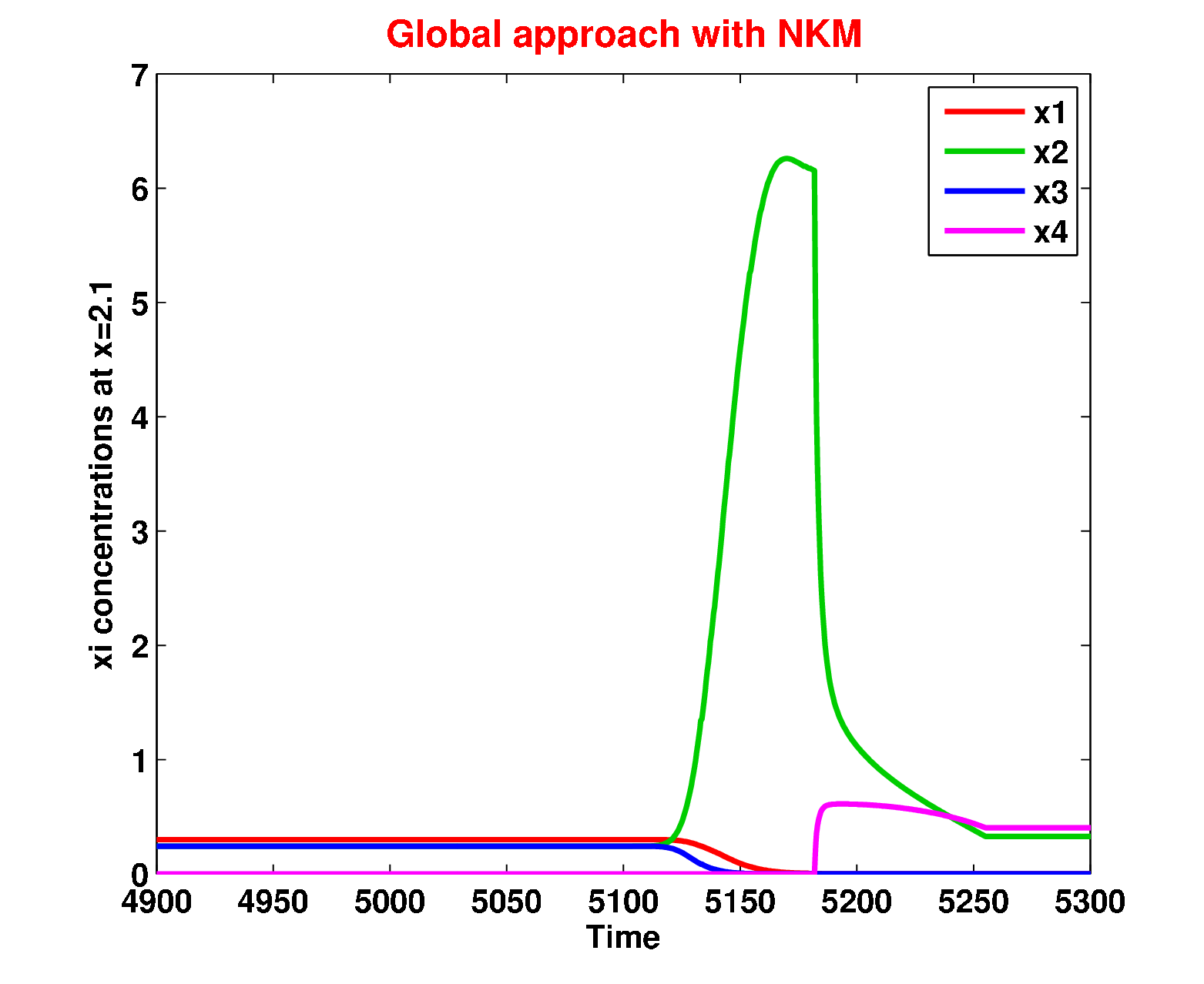}\label{fig:x1x4tpsR}
}
  \caption{Concentration of the components X1 and X4 at the end of the
    column ($x=2.1$) as a function of time}
  \label{fig:x1x4tps}
\end{figure}

First, results showing the evolution of the component species at
various times, and using several spatial and temporal resolutions are
shown on figure~\ref{fig:x1x4L}. The left figure is at time $t=10$,
the right one at $t=50$. As expected, the concentrations remain almost
constant in the middle (reactive) region. Meshes with $220$, $440$,
$660$ and $880$ points have been used, and in each case the time step
is chosen as 0.9 times the limit fixed by the CFL condition. For these
early times, the dependence on the mesh is not very strong.
Elution curves (concentrations at the end of the column as functions
of time) are shown on figure~\ref{fig:x1x4tps}, first for $t$ going from 0 to 400
(figure~\ref{fig:x1x4tpsL}), then for $t$ going from 4900 to 5300
(figure~\ref{fig:x1x4tpsR}). 
The elution curves show that the correct limiting behavior is reached
before the leaching phase begins.

\medskip
The output results required in the benchmark definition are included.
Most were obtained with a 220 points mesh, which may not be sufficient,
as will be seen below. It has not yet been possible to obtain results
with a finer mesh resolution for significantly longer times. 

\begin{figure}[ht]
  \centering
  \subfloat[Total dissolved concentration C3]%
  {\includegraphics[width=0.31\textwidth]{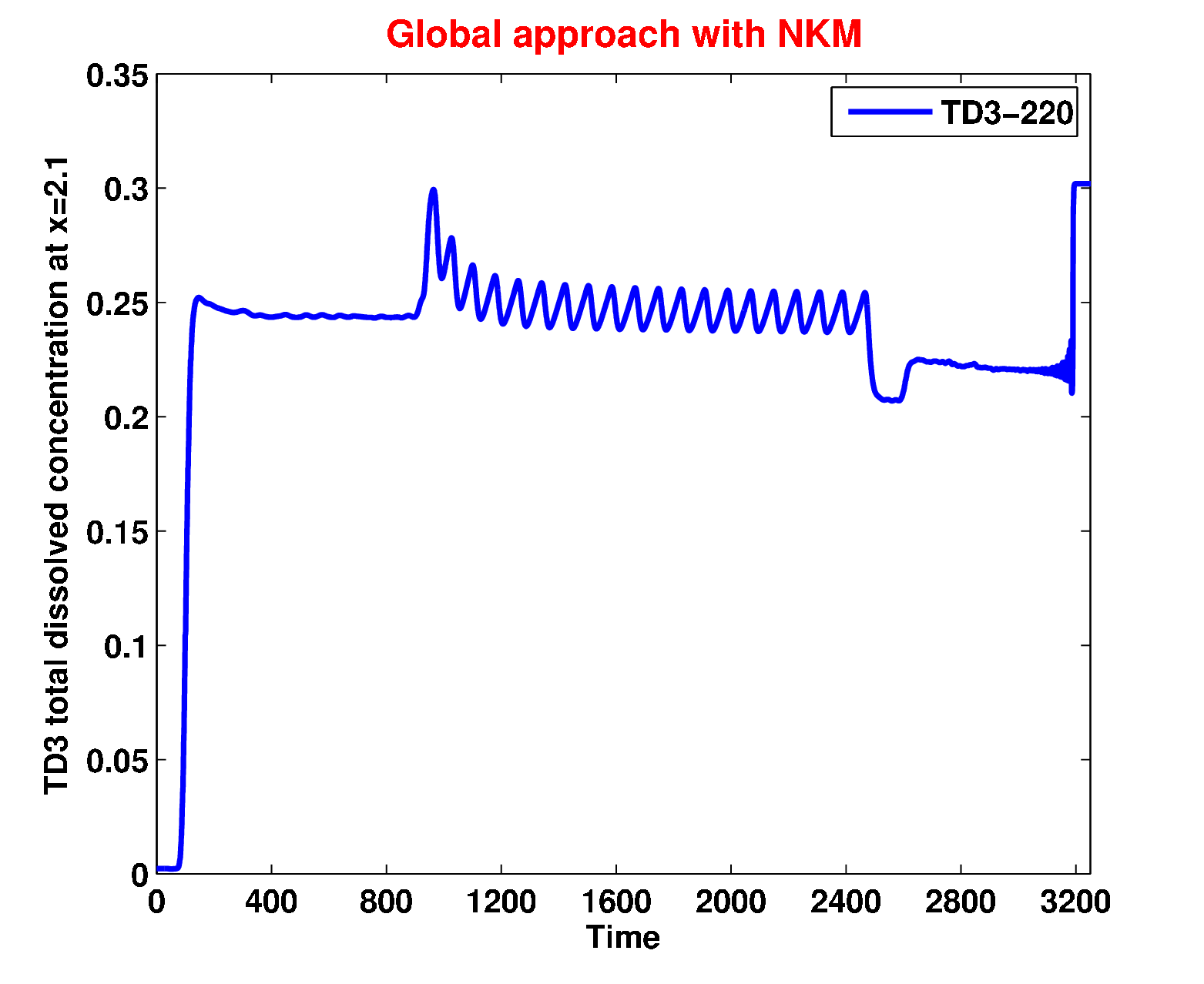}\label{fig:fig123L}}
  \subfloat[Component X3]%
  {\includegraphics[width=0.31\textwidth]{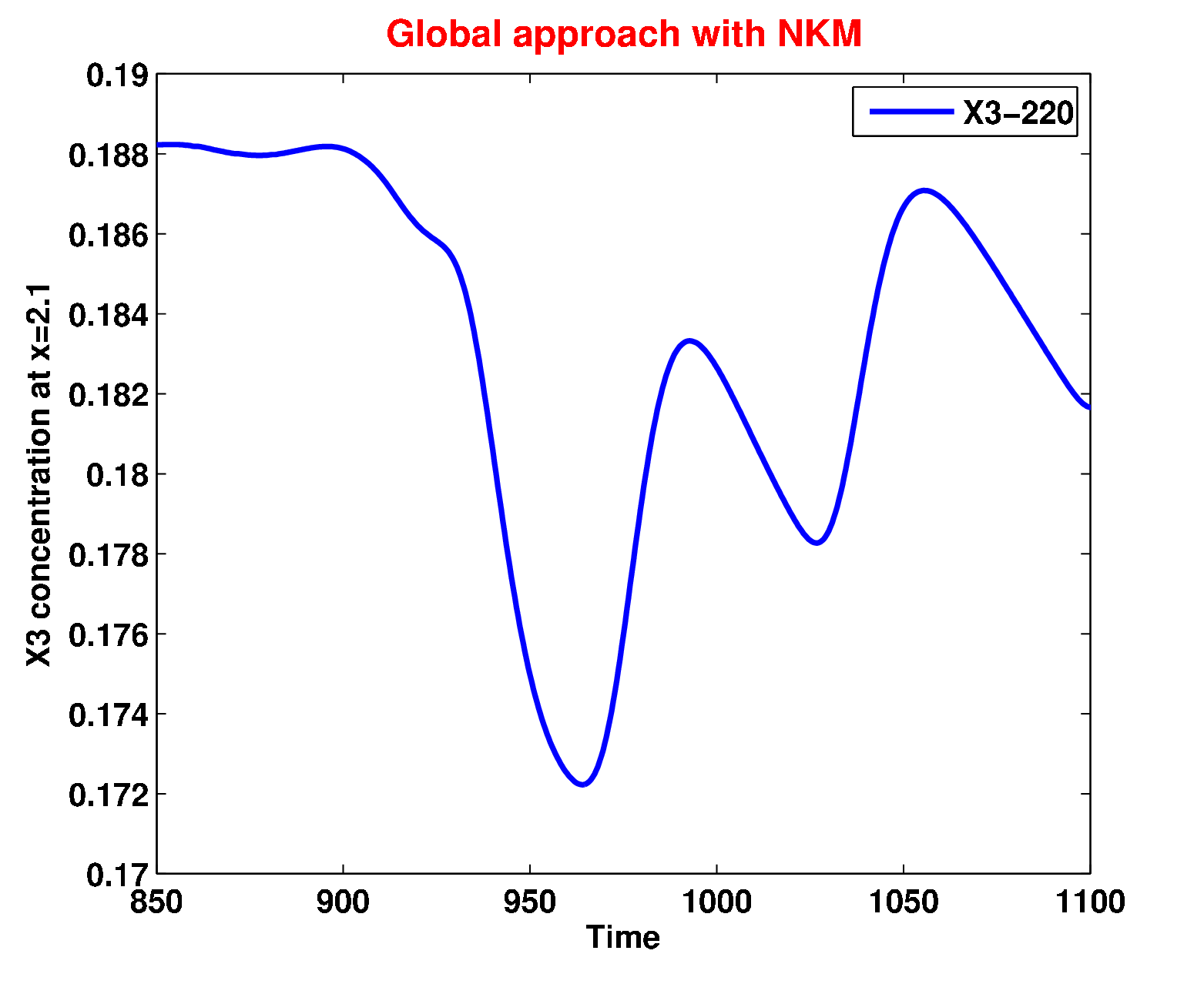}\label{fig:fig123M}}
  \subfloat[Species C1]%
  {\includegraphics[width=0.31\textwidth]{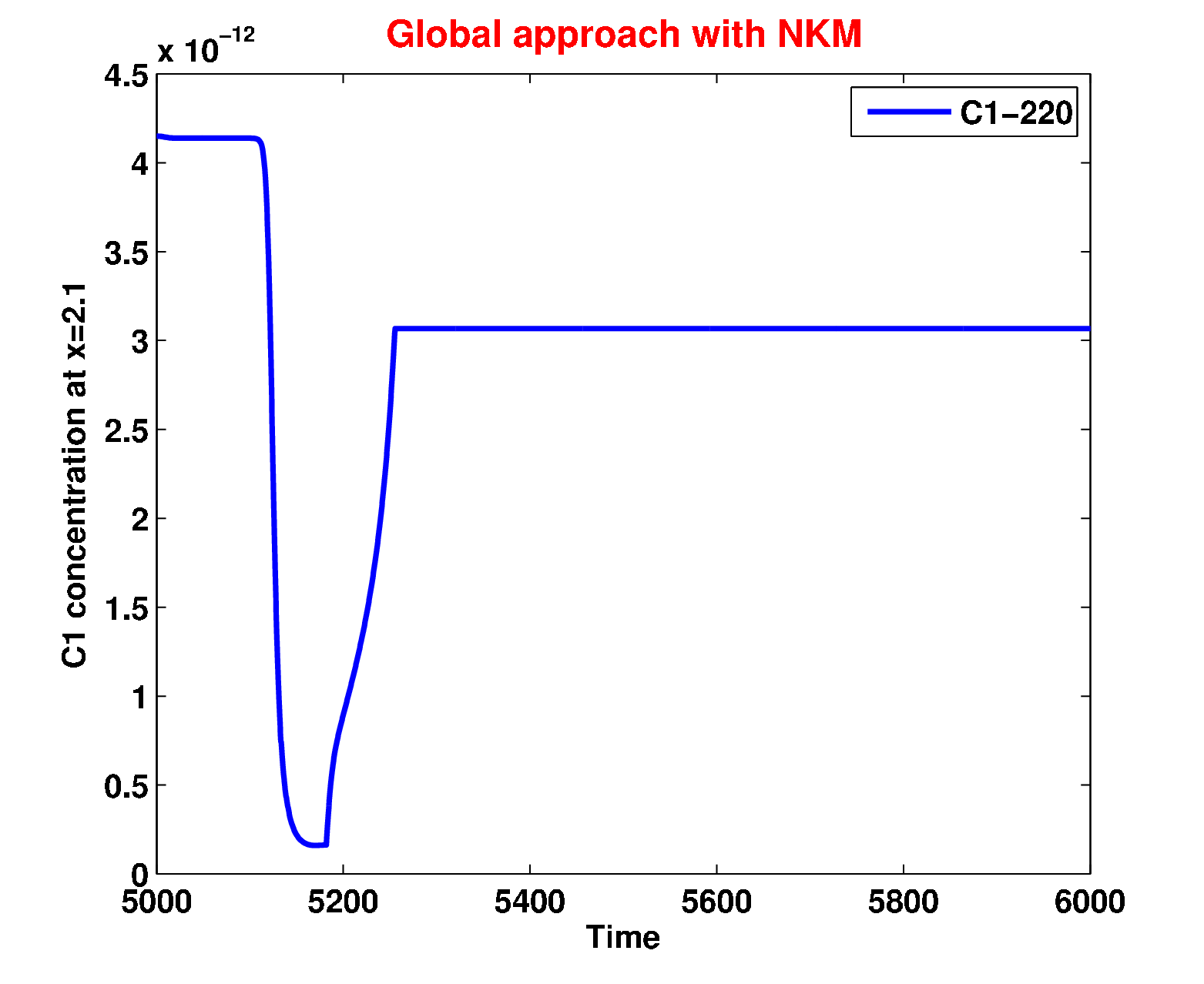}\label{fig:fig123R}}
  \caption{Elution curve (concentrations at $x =  2.1$ as as function
    of time)}
  \label{fig:fig123}
\end{figure}

Figures~\ref{fig:fig123L}, figures~\ref{fig:fig123M}
and~\ref{fig:fig123R} (elution curve for the
total dissolved concentration of component X3,a nd species C1) show an
oscillations pattern that has been observed by other groups working on
the benchmark. These oscillations have been convincingly explained by
V. Lagneau~\cite{lagvdl:09} as being due to the interaction of the very rapid chemistry
and the discrete nature of the grid. They are a discretization artifact,
but appear independently of the method. They can be reduced by using a
more refined grid.

\begin{figure}[ht]
  \centering
  \subfloat[Species X1, $t=10$]%
  {\includegraphics[width=0.49\textwidth]{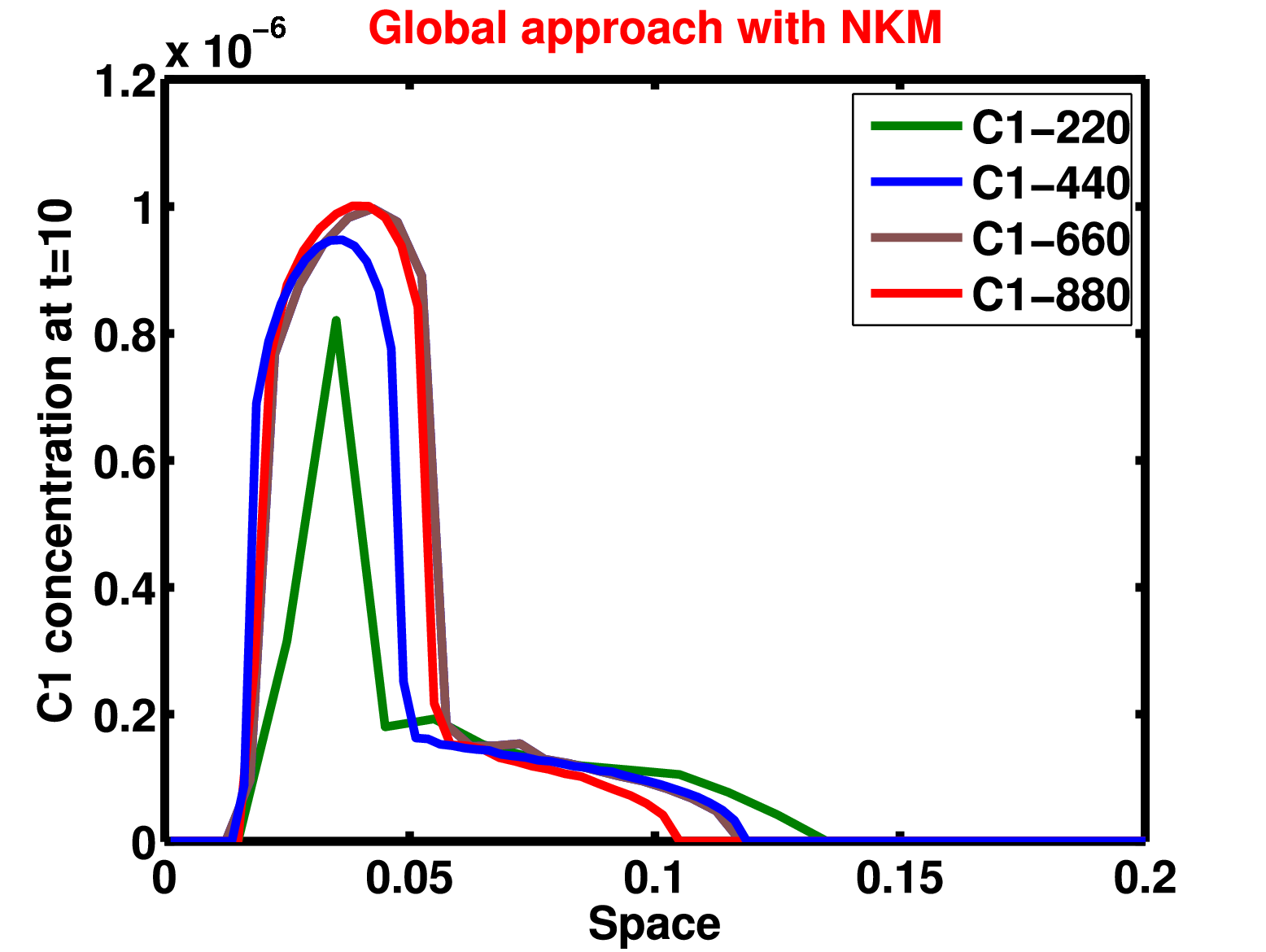}\label{fig:fig45L}}
  \subfloat[Component S, $t=10$]%
  {\includegraphics[width=0.51\textwidth]{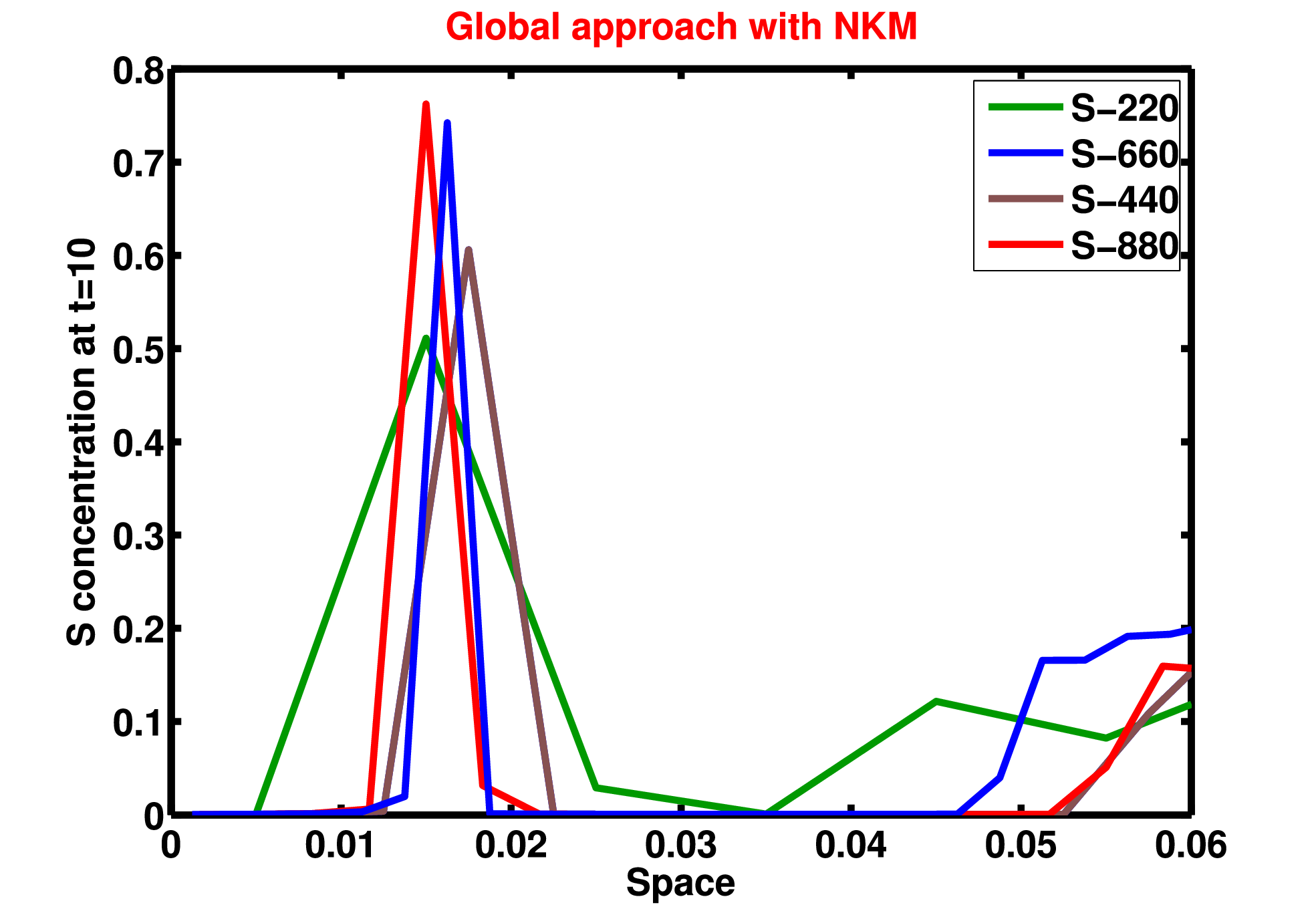}\label{fig:fig45M}}
  \caption{Concentration profiles}
  \label{fig:fig45}
\end{figure}

Figures~\ref{fig:fig45L} and~\ref{fig:fig45M} show the influence on
the mesh, by showing the concentration over a small spatial region,
for time $t=10$ . The concentrations are computed with 4 meshes of increasing
resolution. The peaks in the solution are not resolved satisfactorily
for the coarser mesh, with 220 points, but 660 (and better 880) points
give the correct location and amplitude. Even if the method as it is currently
implemented cannot yet be considered as robust, its ability to locate
these solution features with reasonably coarse meshes was seen as one
of its strong points.  Unfortunately, this may still not be enough to
eliminate the oscillations shown on figure~\ref{fig:fig123}. This
issue is currently being worked on, part of the difficulty being that
increasing the mesh resolution may not be sufficient. As the nonlinear
problem becomes more difficult, it may be necessary to increase the
maximum number of iterations allowed to make sure the Newton--Krylov
method has converged.

% Last, figure~\ref{fig:fig6} shows the concentration of species $X2$ at
% time $t=5010$. 
% \begin{figure}[ht]
%   \centering
%   {\includegraphics[width=0.40\textwidth]{figure6_GL}\label{fig:fig456R}}
%   \caption{Species X2, $t=5010$}
%   \label{fig:fig6}
% \end{figure}

\subsubsection{Performance of the method}
\label{sec:perfs}
The benchmark was intended to be a difficult test for numerical
methods, and this is indeed the case. On the average, more than 20 Newton
iterations are required at each time step, and between 15 and 40
conjugate gradient steps are needed at each nonlinear iterations.

\begin{figure}[bh]
  \centering
  \includegraphics[width=0.5\textwidth]{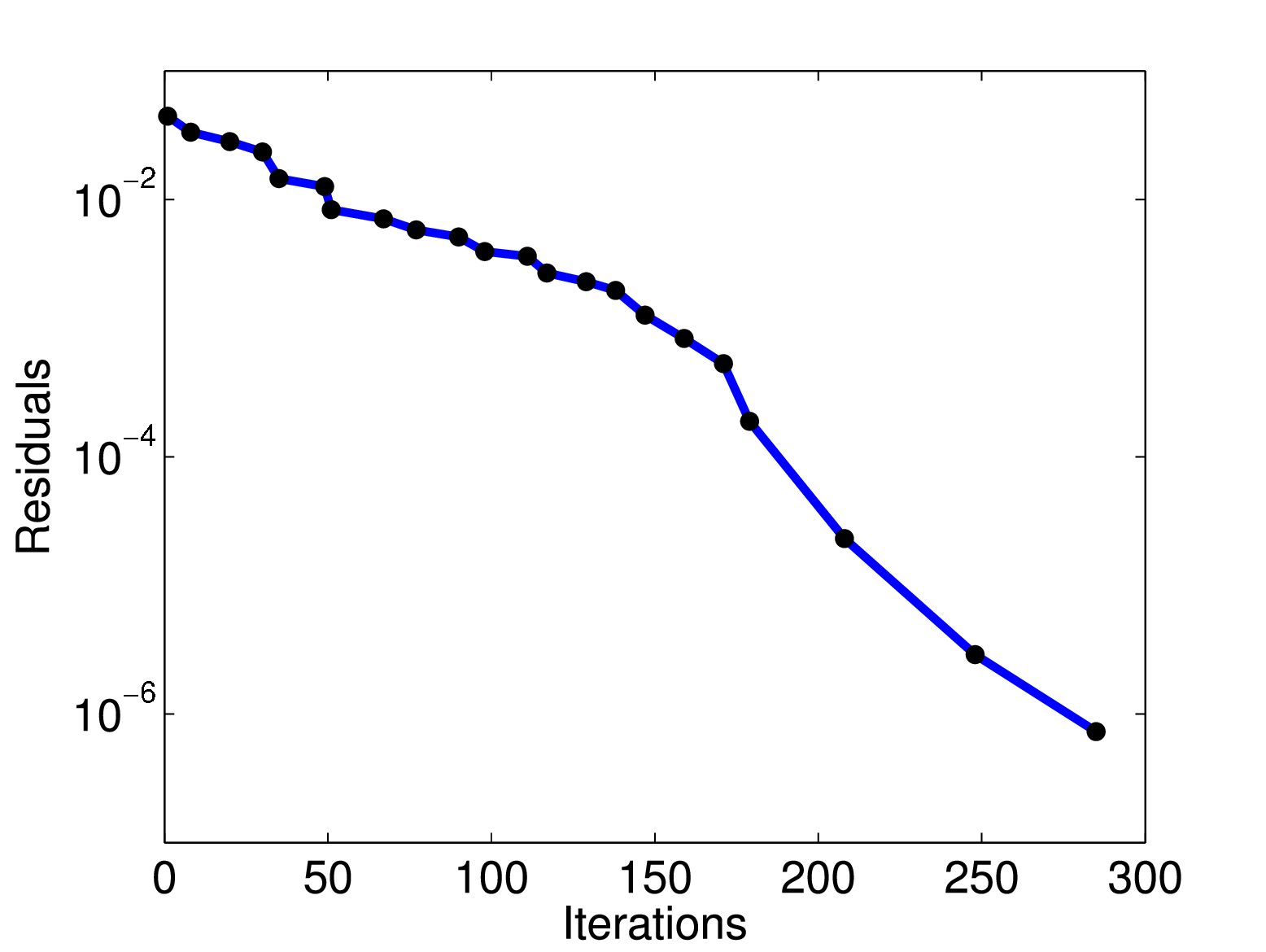}
  \caption{Iterations}
  \label{fig:iters_benchNL}
\end{figure}

Figure~\ref{fig:iters_benchNL} shows a typical time step: the solid
curve shows the cumulative number of conjugate gradient
(alternatively, the number of matrix vector products), and the dots
represent the nonlinear iterations.

Statistics for a single time step are gathered in
table~\ref{tab:statsNK}, for three different mesh resolutions (220,
440 and 660 points). They give the number of non-linear iterations
(NNI) for a (typical) time step, and the total number of linear
iterations (NLI) accummulated over the whole Newton iteration. The
number of nonlinear iterations depends only weakly on the mesh
resolution, whereas the number if linear iterations increases with the
mesh resolution.
\begin{table}[H]
\begin{center}
\begin{tabular}{|c|c|c|c|c|c|}
\hline
\multicolumn{2}{|c|}{\bf Mesh 220}& \multicolumn{2}{c|}{\bf Mesh 440}& \multicolumn{2}{c|}
{\bf Mesh 660}\\        
\hline
NNI& NLI & NNI & NLI & NNI & NLI\\
25   &494   &18   &551      &25   &636 \\ \hline
\end{tabular}
\end{center}
\caption{Statistics on Newton and GNRES iterations, for one time
  step (NNI= Number of NonLinear Iterations, NLI= Number of Linear Iterations).}
\label{tab:statsNK}
\end{table}

Table~\ref{tab:statsNK} shows that the solver spends a large
proportion of its time in the linear solver, despite the adaptive
choice of the forcing parameter (equation~\eqref{eq:forcing}).
Moreover, the number of linear iterations for each nonlinear iteration
also increases with the mesh resolution. Actually, this is expected,
as the solution of the linearized problem includes the solution of the
transport operator, wich has an elliptic-like structure, so that its
condition number grows like the square of the number of grid points.
This problem could be alleviated by using a suitable preconditioner
that would make the number of iterations independent of the mesh
resolution (a domain decomposition preconditioner could be used as
in~\cite{Achdou2000145}). As noticed by Hammond et
al.~\cite{hammvallicht05}, designed a matrix-free preconditioner (so
as to be compatible with the Newton-Krylov framework) is a challenge.
Natural choices would exploit the block structure of the Jacobian, the
simpler ones being based on block-Jacobi, or block Gauss-Seidel.
Operator-splitting as a preconditioner has also been proposed
in~\cite{hammvallicht05}. These possibilities are currently being
explored, exploiting the block structure of the Jacobian, and the
results will be reported in a forthcoming paper~\cite{taakkern:09}.

\section{Conclusions-- Perspectives}

In this paper, it was shown that a global method for coupling
transport with chemistry based on the Newton-Krylov technology can be
implemented while keeping the transport and chemical solvers
separated. The results shown are promising: it is possible to solve
efficiently geochemical problems using the method, although there
remains several issues that need to be addressed.
\begin{itemize}
\item The first is to run test cases on more demanding configurations,
  where the method can be expected to show its full potential. This
  includes the other MoMaS test cases, with a more complex chemistry
  model, and also an implementation of the method in 2 and 3
  dimensions. 
\item It will then certainly be  necessary to explore the question of
  how to precondition the Jacobian, in order to reduce the number of Krylov
  iterations. An natural avenue is to reuse the operator splitting
  methods, as proposed by~\cite{hammvallicht05}. A similar study is
  being carried out for a related, but simpler model,
  see~\cite{taakkern:09}. 
\item The results reported above used a fixed time step, which was
  clearly insufficient for the large interval of integration. To
  successfully solve difficult problems like the benchmark above, it
  will clearly be necessary to use adaptive time stepping. 
\item A more difficult problem will be to take into account
  precipitation--dissolution phenomena in the chemical model. As the
  models are non-differentiable, this makes it more difficult to
  employ Newton's method.
\end{itemize}

As was apparent from the numerical experiments, the method also shows
some limitations. The most serious is its high cost, as each
evaluation of the residual involves the solution of a chemical problem
at each grid point. The fact that the method has two levels of
nonlinear iterations means that it may not be as robust as other
global methods based on a single level of iterations. Finding a good
preconditioner may not be a limitation, but most strategies will
involve solving more transport problems, which will also incur a high
cost.

\section*{Acknowledgments} 
The first author would like to express her
sincere thanks to ITASCA Consultants, France for providing the
necessary Ph.D fellowship to carry out this research in
INRIA-Rocquencourt, France. The second author's work was supported by
Groupement MoMaS CNRS-2439. We gratefully acknowledge sponsorship of GDR
MoMAS by ANDRA, BRGM, CEA, EDF and IRSN. Both authors thank the
referees for their detailed comments, which led to significant
improvements in the contents of the paper.

\bibliography{AmirKern}
\bibliographystyle{spbasic}

\end{document}